# Existence and Smoothness of Navier-Stokes Equations

Svetlin G. Georgiev [a] and Gal Davidi [b]

[a] *Sorbonne University, Paris, France*
e-mail: svetlingeorgiev1@gmail.com

[b] *CFD FEM Engineering, Kiryat Motzkin, Israel*
e-mail: gal.davidi11@gmail.com

In this paper we propose new method for proving of existence of global solutions for 3D Navier-Stokes equations. The proposed method can be applied for investigation of global solutions for other classes of PDEs.

## 1 Introduction

In this paper we investigate a class of IVP for a class of Navier-Stokes equations for existence of global classical solutions. More precisely, we will study the following Navier-Stokes equations

$$u_t + uu_x + vu_y + wu_z + \frac{1}{\rho} p_x - \nu u_{xx} - \nu u_{yy} - \nu u_{zz} = 0$$

$$v_t + uv_x + vv_y + wv_z + \frac{1}{\rho} p_y - \nu v_{xx} - \nu v_{yy} - \nu v_{zz} = 0 \qquad (1.1)$$

$$w_t + uw_x + vw_y + ww_z + \frac{1}{\rho} p_z - \nu w_{xx} - \nu w_{yy} - \nu w_{zz} = 0$$

$$u_x + v_y + w_z = 0 \qquad in \qquad (0,\infty) \times \mathbb{R}^3,$$

$$u(0,x,y,z) = u_0(x,y,z), \qquad v(0,x,y,z) = v_0(x,y,z),$$

$$w(0,x,y,z) = w_0(x,y,z), \qquad (x,y,z) \in \mathbb{R}^3,$$

where $p,u,v,w\colon [0,\infty) \times \mathbb{R}^3 \to \mathbb{R}$ are unknown, $u_0, v_0, w_0 \in C^\infty(\mathbb{R}^3)$ are given functions. This is a system of partial differential equations that governs the flow of a viscous incompressible fluid. Here $\rho$ is the density, $u$ the velocity vector, $p$ is the pressure. The first three equations of (1.1) are Cauchy's momentum equations where the first term is the accelerating time varying term, the second and third are the convective and the hydrostatic terms respectively. The physical example of the convective term can be described as a river that is converging, the case where the term is increasing and the river diverging the case where the term is decreasing. The hydrostatic term describes flow from high pressure to low pressure. The forth term is the viscosity term with the coefficient $\nu$ the kinematical viscosity. This term describes the ability of the fluid to induce motion of neighboring particles. On the right hand side we have the external forces density term. This term can include: gravity, magneto-hydrodynamic force, and so on. The fourth equation of (1.1) is the nullification of the divergence due to incompressibility condition. Turbulent fluid motions are believed to be well modeled by the Navier-Stokes equations. In the case of the 3D version of the NS equations the existence problem is an unsolved issue([1]).

We recall that global existence of weak solutions of the Navier-Stokes equations is known to hold in every space dimension. Uniqueness of weak solutions and global existence of strong solutions is known in dimension two [4]. In dimension three, global existence of strong solutions of the Navier-Stokes equations in thin three-dimensional domains began with the papers [5] and [6], where is used the methods in [2] and [3].

In this paper we propose new method for investigation of equations (1.1). The proposed method gives existence of classical solutions for the problem (1.1).

Without loss of generality we assume that $\rho = \nu = 1$. For a set $A$ such that $A \subset \mathbb{R}$ or $A \subset \mathbb{R}^2$ or $A \subset \mathbb{R}^3$, with $\mu(A)$ we will denote its measure.

Let $\mathbb{R}^3 = \cup_{j=1}^{\infty} D_j$, where $D_j$ be bounded subsets of $\mathbb{R}^3$ satisfying the following conditions

1. $D_i \cap D_j = \emptyset$ for $i \neq j$,, $i, j \in \{1,2,\dots\}$,
2. $D_j$ and $D_{j+1}$, $j \in \{1,2,\dots\}$, are adjoining,
3. if $(x_1, y_1, z_1) \in D_j$, for some $j \in \{1,2,\dots\}$, is fixed and

$$D_{jx_1} = \{(y,z) \in \mathbb{R}^2 : (x_1, y, z) \in D_j\},$$

$$D_{jy_1} = \{(x,z) \in \mathbb{R}^2 : (x, y_1, z) \in D_j\},$$

$$D_{jz_1} = \{(x,y) \in \mathbb{R}^2 : (x, y, z_1) \in \mathrm{D}_j\},$$

$$D_{jx_1y_1} = \{z \in \mathbb{R} : (x_1, y_1, z) \in D_j\},$$

$$D_{jx_1z_1} = \{y \in \mathbb{R} : (x_1, y, z_1) \in D_j\},$$

$$D_{jy_1z_1} = \{x \in \mathbb{R} : (x, y_1, z_1) \in D_j\},$$

then

$$\mu(D_{jx_1}) \leq \mu(D_j),$$

$$\mu(D_{jy_1}) \leq \mu(D_j),$$

$$\mu(D_{jz_1}) \leq \mu(D_j),$$

$$\mu(D_{jx_1y_1}) \leq \mu(D_j),$$

$$\mu(D_{jx_1z_1}) \leq \mu(D_j),$$

$$\mu(D_{jy_1z_1}) \leq \mu(D_j), \qquad j \in \{1,2,\dots\}.$$

4. For any $j \in \mathbb{N}$, let $D_{jj}$ is a compact subset of $D_j$ such that $D_{jj} \neq D_j$ and there exists $(x_0, y_0, z_0) \in \partial D_{jj}$ so that if $(x_1, y_1, z_1) \not\subset D_{jj}$, then $[x_0, x_1] \times [y_0, y_1] \times [z_0, z_1] \cap D_{jj} = \emptyset$, where $[x_0, x_1]$ is the segment with end points $x_0$ and $x_1$, $[y_0, y_1]$ is the segment with end points $y_0$ and $y_1$, $[z_0, z_1]$ is the segment with end points $z_0$ and $z_1$.

For a set $B \subset \mathbb{R}^3$ and a function $f: \mathbb{R}^3 \to \mathbb{R}$ with $f|_B$ we denote the restriction of f to $B$.

Our main result is as follows.

**Theorem 1.1** *Let* $u_0, v_0, w_0 \in C^\infty(\mathbb{R}^3)$ *be such that*

$$u_0|_{D_j}, v_0|_{D_j}, w_0|_{D_j} \in C_0^\infty(D_j),$$

$$supp\left(u_0|_{D_j}\right), \qquad supp\left(v_0|_{D_j}\right), \qquad supp\left(w_0|_{D_j}\right) \subset D_{jj}, \qquad j \in \{1,2,\dots\},$$

and

$$|\,\partial_x^{\alpha_1}\,\partial_y^{\alpha_2}\,\partial_z^{\alpha_3} u_0(x,y,z)| \leq C_{\alpha_1\alpha_2\alpha_3K}\left(1+\sqrt{x^2+y^2+z^2}\right)^{-K},$$

$$|\,\partial_x^{\alpha_1}\,\partial_y^{\alpha_2}\,\partial_z^{\alpha_3} v_0(x,y,z)| \leq C_{\alpha_1\alpha_2\alpha_3K}\left(1+\sqrt{x^2+y^2+z^2}\right)^{-K}, \tag{1.2}$$

$$|\,\partial_x^{\alpha_1}\,\partial_y^{\alpha_2}\,\partial_z^{\alpha_3} w_0(x,y,z)| \leq C_{\alpha_1\alpha_2\alpha_3K}\left(1+\sqrt{x^2+y^2+z^2}\right)^{-K},$$

on $\mathbb{R}^3$, for any $\alpha_1,\ \alpha_2,\ \alpha_3 \in \mathbb{N}_0$ and positive constant $K$,

$$\| u_0 \|_{L^2(\mathbb{R}^3)} \leq 1, \qquad \| v_0 \|_{L^2(\mathbb{R}^3)} \leq 1, \qquad \| w_0 \|_{L^2(\mathbb{R}^3)} \leq 1,$$

$$\sup_{(x,y,z)\in D_j} |u_0(x,y,z)| < \frac{1}{2^{j^2}\sqrt{\mu(D_j)}},$$

$$\sup_{(x,y,z)\in D_j} |u_{0x}(x,y,z)| < \frac{1}{2^{j^2}\sqrt{\mu(D_j)}},$$

$$\sup_{(x,y,z)\in D_j} |u_{0xx}(x,y,z)| < \frac{1}{2^{j^2}\sqrt{\mu(D_j)}},$$

$$\sup_{(x,y,z)\in D_j} |u_{0y}(x,y,z)| < \frac{1}{2^{j^2}\sqrt{\mu(D_j)}}, \sup_{(x,y,z)\in D_j} |u_{0yy}(x,y,z)| < \frac{1}{2^{j^2}\sqrt{\mu(D_j)}},$$

$$\sup_{(x,y,z)\in D_j} |u_{0z}(x,y,z)| < \frac{1}{2^{j^2}\sqrt{\mu(D_j)}}, \qquad \sup_{(x,y,z)\in D_j} |u_{0zz}(t,x,y,z)| < \frac{1}{2^{j^2}\sqrt{\mu(D_j)}},$$

$$\sup_{(x,y,z)\in D_j} |v_0(x,y,z)| < \frac{1}{2^{j^2}\sqrt{\mu(D_j)}},$$

$$\sup_{(x,y,z)\in D_j} |v_{0x}(x,y,z)| < \frac{1}{2^{j^2}\sqrt{\mu(D_j)}},$$

$$\sup_{(x,y,z)\in D_j} |v_{0xx}(x,y,z)| < \frac{1}{2^{j^2}\sqrt{\mu(D_j)}},$$

$$\sup_{(x,y,z)\in D_j} |v_{0y}(x,y,z)| < \frac{1}{2^{j^2}\sqrt{\mu(D_j)}}, \sup_{(x,y,z)\in D_j} |v_{0yy}(x,y,z)| < \frac{1}{2^{j^2}\sqrt{\mu(D_j)}},$$

$$\sup_{(x,y,z)\in D_j} |v_{0z}(x,y,z)| < \frac{1}{2^{j^2}\sqrt{\mu(D_j)}}, \qquad \sup_{(x,y,z)\in D_j} |v_{0zz}(t,x,y,z)| < \frac{1}{2^{j^2}\sqrt{\mu(D_j)}},$$

$$\sup_{(x,y,z)\in D_j} |w_0(x,y,z)| < \frac{1}{2^{j^2}\sqrt{\mu(D_j)}},$$

$$\sup_{(x,y,z)\in D_j} |w_{0x}(x,y,z)| < \frac{1}{2^{j^2}\sqrt{\mu(D_j)}},$$

$$\sup_{(x,y,z)\in D_j} |w_{0xx}(x,y,z)| < \frac{1}{2^{j^2}\sqrt{\mu(D_j)}},$$

$$\sup_{(x,y,z)\in D_j} |w_{0y}(x,y,z)| < \frac{1}{2^{j^2}\sqrt{\mu(D_j)}}, \sup_{(x,y,z)\in D_j} |w_{0yy}(x,y,z)| < \frac{1}{2^{j^2}\sqrt{\mu(D_j)}},$$

$$\sup_{(x,y,z)\in D_j} |w_{0z}(x,y,z)| < \frac{1}{2^{j^2}\sqrt{\mu(D_j)}}, \qquad \sup_{(x,y,z)\in D_j} |w_{0zz}(t,x,y,z)| < \frac{1}{2^{j^2}\sqrt{\mu(D_j)}},$$

$j \in \mathbb{N}$. Then the problem (1.1) has a solution $(u,v,w,p) \in (C^\infty([0,\infty)\times\mathbb{R}^3))^4$ such that

$$\int_{\mathbb{R}^3} |u(t,x,y,z)|^2 dxdydz \le C_1, \qquad \int_{\mathbb{R}^3} |v(t,x,y,z)|^2 dxdydz \le C_1,$$

$$\int_{\mathbb{R}^3} |w(t,x,y,z)|^2 dxdydz \le C_1, \qquad \int_{\mathbb{R}^3} |p(t,x,y,z)|^2 dxdydz \le C_1,$$

for some constant $C_1 > 1$ and for any $t \in [0,\infty)$.

**Remark 1.2** *If $u_0 \not\equiv 0$, $v_0 \not\equiv 0$, $w_0 \not\equiv 0$, then we obtain a nontrivial solution of the system (1.1).*

The paper is organized as follows. In the next section, we give some preliminary results which we will use for the proof of our main result. In Section 3 we prove our main result.

## 2 Preliminaries

**Definition 2.1** *Let $(X,d)$ be a metric space and $M$ be a subset of $X$. The mapping $T: M \to X$ is said to be expansive if there exists a constant $h > 1$ such that*

$$d(Tx,Ty) \ge hd(x,y)$$

for any $x, y \in M$.

**Theorem 2.2 (, Theorem 2.4)** *Let $X$ be a nonempty closed convex subset of a Banach space $E$. Suppose that $T$ and $S$ map $X$ into $E$ such that*

1. $S$ is continuous and $S(X)$ resides in a compact subset of $E$.
2. $T: X \to E$ is expansive.
3. $S(X) \subset (I - T)(E)$ and $[x = Tx + Sy, y \in X] \Rightarrow x \in X$(or $S(X) \subset (I - T)(X)$).

Then there exists a point $x^* \in X$ such that

$$Sx^* + Tx^* = x^*.$$

The following result is a consequence of Theorem 2.2.

**Theorem 2.3** *Let $X$ be a nonempty closed convex subset of a Banach space $E$ and $Y$ is a nonempty compact subset of $E$ such that $X \subset Y$, $Y \neq X$. Suppose that $T$ and $S$ map $X$ into $E$ such that*

1. $S$ is continuous and $S(X)$ resides in $Y$.
2. $T: X \to E$ is linear, continuous and expansive, and $T: X \to Y$ is onto, and $\{x - z : x \in X, \quad z \in S(X)\} \subset Y$.

Then there exists an $x^* \in X$ such that

$$Tx^* + Sx^* = x^*.$$

*Proof.* Since $Y$ is compact and $S(X)$ resides in $Y$, we have that the first condition of Theorem 2.2 holds. Because $T: X \to E$ is expansive, we have that the second condition of Theorem 2.2 holds. Note that $T^{-1}: Y \to E$ exists, it is linear and contractive with a constant $l \in (0,1)$. Let $z \in S(X)$ be fixed. Take $y_0 \in Y$ arbitrarily. Define the sequence $\{y_n\}_{n \in \mathbb{N}}$ as follows.

$$y_{n+1} = T^{-1} y_n - z, \qquad n \in \mathbb{N} \cup \{0\}.$$

Note that $\{y_n\}_{n \in \mathbb{N}} \subset Y$ because $T^{-1} y_n \in X$ and $z \in S(X)$. Then

$$\| y_2 - y_1 \| = \| T^{-1} y_1 - T^{-1} y_0 \|$$

$$\leq l \| y_1 - y_0 \|,$$

$$\| y_3 - y_2 \| = \| T^{-1} y_2 - T^{-1} y_1 \|$$

$$\leq l \| y_2 - y_1 \|$$

$$\leq l^2 \| y_1 - y_0 \|.$$

Using the principle of the mathematical induction, we get

$$\| y_{n+1} - y_n \| \leq l^n \| y_1 - y_0 \|, \qquad n \in \mathbb{N},$$

where $\mathbb{N}$ is the set of the natural numbers. Now, for $m > n$, $m, n \in \mathbb{N}$, we find

$$\| y_m - y_n \| \leq \| y_m - y_{m-1} \| + \cdots + \| y_{n+1} - y_n \|$$

$$\leq (l^{m-1} + \cdots + l^n) \| y_1 - y_0 \|$$

$$\leq l^n \sum_{j=0}^{\infty} l^j \parallel y_1 - y_0 \parallel$$

$$= \frac{l^n}{1-l} \parallel y_1 - y_0 \parallel.$$

Therefore $\{y_n\}_{n\in\mathbb{N}}$ is a Cauchy sequence of elements of $Y \subset E$. Since $E$ is a Banach space, it follows that the sequence $\{y_n\}_{n\in\mathbb{N}}$ is convergent to an element $y^* \in E$. Because $\{y_n\}_{n\in\mathbb{N}} \subset Y$ and $Y \subset E$ is compact, we have that $y^* \in Y$. Thus

$$y^* = T^{-1}y^* - z$$

or

$$z^* = Tz^* + z, \qquad z^* = T^{-1}y^* \in X.$$

Because $z \in S(X)$ was arbitrarily chosen, we conclude that $S(X) \subset (I - T)(X)$, i.e., the third condition of Theorem 2.2 holds. Hence Theorem 2.2, it follows that there exists an $x^* \in X$ such that

$$Tx^* + Sx^* = x^*.$$

This completes the proof.

## 3 Proof of the Main Result

From the fourth equation of the system (1.1), we get

$$u(u_x + v_y + w_z) = 0, \qquad v(u_x + v_y + w_z) = 0, \qquad w(u_x + v_y + w_z) = 0.$$

Then the system (1.1) we can rewrite in the form

$$u_t + uu_x + vu_y + wu_z + u(u_x + v_y + w_z) + p_x - u_{xx} - u_{yy} - u_{zz} = 0$$

$$v_t + uv_x + vv_y + wv_z + v(u_x + v_y + w_z) + p_y - v_{xx} - v_{yy} - v_{zz} = 0$$

$$w_t + uw_x + vw_y + ww_z + w(u_x + v_y + w_z) + p_z - w_{xx} - w_{yy} - w_{zz} = 0$$

$$u_x + v_y + w_z = 0 \qquad in \qquad (0,\infty) \times \mathbb{R}^3,$$

$$u(0,x,y,z) = u_0(x,y,z), \qquad v(0,x,y,z) = v_0(x,y,z),$$

$$w(0,x,y,z) = w_0(x,y,z), \qquad (x,y,z) \in \mathbb{R}^3,$$

whereupon

$$u_t + (u^2)_x + (uv)_y + (uw)_z + p_x - u_{xx} - u_{yy} - u_{zz} = 0$$

$$v_t + (uv)_x + (v^2)_y + (vw)_z + p_y - v_{xx} - v_{yy} - v_{zz} = 0$$

$$w_t + (uw)_x + (vw)_y + (w^2)_z + p_z - w_{xx} - w_{yy} - w_{zz} = 0$$

$$u_x + v_y + w_z = 0 \qquad in \qquad (0,\infty) \times \mathbb{R}^3, \tag{3.1}$$

$$u(0,x,y,z) = u_0(x,y,z), \qquad v(0,x,y,z) = v_0(x,y,z),$$

$$w(0,x,y,z) = w_0(x,y,z), \qquad (x,y,z) \in \mathbb{R}^3.$$

**Remark 3.1** *We note that using the fourth equation of (3.1) we can obtain the system (1.1).*

Let $j \in \{1,2,\dots\}$ be arbitrarily chosen. Firstly, we consider the problem

$$
\begin{aligned}
& u_t + (u^2)_x + (uv)_y + (uw)_z + p_x - u_{xx} - u_{yy} - u_{zz} = 0 \\
& v_t + (uv)_x + (v^2)_y + (vw)_z + p_y - v_{xx} - v_{yy} - v_{zz} = 0 \\
& w_t + (uw)_x + (vw)_y + (w^2)_z + p_z - w_{xx} - w_{yy} - w_{zz} = 0 \\
& u_x + v_y + w_z = 0 \qquad \text{in} \qquad (0,1] \times D_j, \\
& u(0,x,y,z) = u_0(x,y,z), \qquad v(0,x,y,z) = v_0(x,y,z), \\
& w(0,x,y,z) = w_0(x,y,z), \qquad (x,y,z) \in D_j.
\end{aligned}
\tag{3.2}
$$

We will prove that the problem (3.2) has a solution $(u,v,w,p)$ such that $u,v,w,p \in C^1([0,1], C_0^2(D_j))$.

For $u,v,w,p \in C^1([0,1], C_0^2(D_j))$, we define

$$
\begin{aligned}
I_1^{1j}(u,v,w,p) = & \int_{x_0}^{x} \int_{x_0}^{x_1} \int_{y_0}^{y} \int_{y_0}^{y_1} \int_{z_0}^{z} \int_{z_0}^{z_1} (u(t,\alpha,\beta,\gamma) \\
& - u_0(\alpha,\beta,\gamma)) d\gamma dz_1 d\beta dy_1 d\alpha dx_1 \\
& + \int_0^t \int_{x_0}^{x} \int_{y_0}^{y} \int_{y_0}^{y_1} \int_{z_0}^{z} \int_{z_0}^{z_1} u^2(s,\alpha,\beta,\gamma) d\gamma dz_1 d\beta dy_1 d\alpha ds \\
& + \int_0^t \int_{x_0}^{x} \int_{x_0}^{x_1} \int_{y_0}^{y} \int_{z_0}^{z} \int_{z_0}^{z_1} u(s,\alpha,\beta,\gamma) v(s,\alpha,\beta,\gamma) d\gamma dz_1 d\beta d\alpha dx_1 ds \\
& + \int_0^t \int_{x_0}^{x} \int_{x_0}^{x_1} \int_{y_0}^{y} \int_{y_0}^{y_1} \int_{z_0}^{z} u(s,\alpha,\beta,\gamma) w(s,\alpha,\beta,\gamma) d\gamma d\beta dy_1 d\alpha dx_1 ds \\
& + \int_0^t \int_{x_0}^{x} \int_{y_0}^{y} \int_{y_0}^{y_1} \int_{z_0}^{z} \int_{z_0}^{z_1} p(s,\alpha,\beta,\gamma) d\gamma dz_1 d\beta dy_1 d\alpha ds \\
& - \int_0^t \int_{y_0}^{y} \int_{y_0}^{y_1} \int_{z_0}^{z} \int_{z_0}^{z_1} u(s,x,\beta,\gamma) d\gamma dz_1 d\beta dy_1 ds \\
& - \int_0^t \int_{x_0}^{x} \int_{x_0}^{x_1} \int_{z_0}^{z} \int_{z_0}^{z_1} u(s,\alpha,y,\gamma) d\gamma dz_1 d\alpha dx_1 ds \\
& - \int_0^t \int_{x_0}^{x} \int_{x_0}^{x_1} \int_{y_0}^{y} \int_{y_0}^{y_1} u(s,\alpha,\beta,z) d\beta dy_1 d\alpha dx_1 ds,
\end{aligned}
$$

$$I_2^{1j}(u,v,w,p) = \int_{x_0}^{x}\int_{x_0}^{x_1}\int_{y_0}^{y}\int_{y_0}^{y_1}\int_{z_0}^{z}\int_{z_0}^{z_1}(v(t,\alpha,\beta,\gamma)$$

$$-v_0(\alpha,\beta,\gamma))d\gamma dz_1 d\beta dy_1 d\alpha dx_1$$

$$+\int_0^t\int_{x_0}^{x}\int_{y_0}^{y}\int_{y_0}^{y_1}\int_{z_0}^{z}\int_{z_0}^{z_1} u(s,\alpha,\beta,\gamma)v(s,\alpha,\beta,\gamma)d\gamma dz_1 d\beta dy_1 d\alpha ds$$

$$+\int_0^t\int_{x_0}^{x}\int_{x_0}^{x_1}\int_{y_0}^{y}\int_{z_0}^{z}\int_{z_0}^{z_1} v^2(s,\alpha,\beta,\gamma)d\gamma dz_1 d\beta d\alpha dx_1 ds$$

$$+\int_0^t\int_{x_0}^{x}\int_{x_0}^{x_1}\int_{y_0}^{y}\int_{y_0}^{y_1}\int_{z_0}^{z} v(s,\alpha,\beta,\gamma)w(s,\alpha,\beta,\gamma)d\gamma d\beta dy_1 d\alpha dx_1 \mathrm{ds}$$

$$+\int_0^t\int_{x_0}^{x}\int_{x_0}^{x_1}\int_{y_0}^{y}\int_{z_0}^{z}\int_{z_0}^{z_1} p(s,\alpha,\beta,\gamma)d\gamma dz_1 d\beta d\alpha dx_1 ds$$

$$-\int_0^t\int_{y_0}^{y}\int_{y_0}^{y_1}\int_{z_0}^{z}\int_{z_0}^{z_1} v(s,x,\beta,\gamma)d\gamma dz_1 d\beta dy_1 ds$$

$$-\int_0^t\int_{x_0}^{x}\int_{x_0}^{x_1}\int_{z_0}^{z}\int_{z_0}^{z_1} v(s,\alpha,y,\gamma)d\gamma dz_1 d\alpha dx_1 ds$$

$$-\int_0^t\int_{x_0}^{x}\int_{x_0}^{x_1}\int_{y_0}^{y}\int_{y_0}^{y_1} v(s,\alpha,\beta,z)d\beta dy_1 d\alpha dx_1 ds,$$

$$I_3^{1j}(u,v,w,p) = \int_{x_0}^{x}\int_{x_0}^{x_1}\int_{y_0}^{y}\int_{y_0}^{y_1}\int_{z_0}^{z}\int_{z_0}^{z_1}(w(t,\alpha,\beta,\gamma)$$

$$-w_0(\alpha,\beta,\gamma))d\gamma dz_1 d\beta dy_1 d\gamma dx_1$$

$$+\int_0^t\int_{x_0}^{x}\int_{y_0}^{y}\int_{y_0}^{y_1}\int_{z_0}^{z}\int_{z_0}^{z_1} u(s,\alpha,\beta,\gamma)w(s,\alpha,\beta,\gamma)d\gamma dz_1 d\beta dy_1 d\alpha ds$$

$$+\int_0^t\int_{x_0}^{x}\int_{x_0}^{x_1}\int_{y_0}^{y}\int_{z_0}^{z}\int_{z_0}^{z_1} v(s,\alpha,\beta,\gamma)w(s,\alpha,\beta,\gamma)d\gamma dz_1 d\beta d\alpha dx_1 ds$$

$$+\int_0^t\int_{x_0}^{x}\int_{x_0}^{x_1}\int_{y_0}^{y}\int_{y_0}^{y_1}\int_{z_0}^{z} w^2(s,\alpha,\beta,\gamma)d\gamma d\beta dy_1 d\alpha dx_1 ds$$

$$+\int_0^t\int_{x_0}^{x}\int_{x_0}^{x_1}\int_{y_0}^{y}\int_{y_0}^{y_1}\int_{z_0}^{z} p(s,\alpha,\beta,\gamma)d\gamma d\beta dy_1 d\alpha dx_1 ds$$

$$-\int_0^t\int_{y_0}^{y}\int_{y_0}^{y_1}\int_{z_0}^{z}\int_{z_0}^{z_1} w(s,x,\beta,\gamma)d\gamma dz_1 d\beta dy_1 ds$$

$$-\int_0^t\int_{x_0}^{x}\int_{x_0}^{x_1}\int_{z_0}^{z}\int_{z_0}^{z_1} w(s,\alpha,y,\gamma)d\gamma dz_1 d\alpha dx_1 ds$$

$$-\int_0^t\int_{x_0}^{x}\int_{x_0}^{x_1}\int_{y_0}^{y}\int_{y_0}^{y_1} w(s,\alpha,\beta,z)d\beta dy_1 d\alpha dx_1 ds,$$

$$I_4^{1j}(u,v,w,p) = \int_0^t\int_{x_0}^{x}\int_{y_0}^{y}\int_{\mathrm{y}_0}^{y_1}\int_{z_0}^{z}\int_{z_0}^{z_1} u(s,\alpha,\beta,\gamma)d\gamma dz_1 d\beta dy_1 d\alpha ds$$

$$+\int_0^t \int_{x_0}^{x} \int_{x_0}^{x_1} \int_{y_0}^{y} \int_{z_0}^{z} \int_{z_0}^{z_1} v(s,\alpha,\beta,\gamma) d\gamma dz_1 d\beta d\alpha dx_1 ds$$

$$+\int_0^t \int_{x_0}^{x} \int_{x_0}^{x_1} \int_{y_0}^{y} \int_{y_0}^{y_1} \int_{z_0}^{z} w(s,\alpha,\beta,\gamma) d\gamma d\beta dy_1 d\alpha dx_1 ds.$$

Note that the operators $I_k^{1j}$ corresponds to $k$th equation of the problem (3.2), $k \in \{1,2,3,4\}$.

**Lemma 3.2** *Every solution* $(u,v,w,p) \in \left(C^1([0,1],C_0^2(D_j))\right)^4$ *of the system*

$$\begin{aligned} I_1^{1j}(u,v,w,p) &= 0 \\ I_2^{1j}(u,v,w,p) &= 0 \\ I_3^{1j}(u,v,w,p) &= 0 \\ I_4^{1j}(u,v,w,p) &= 0 \end{aligned} \tag{3.3}$$

is a solution of the problem (3.2).

*Proof.* Consider the equation

$$I_1^{1j}(u,v,w,p) = 0 \qquad \text{for} \qquad (t,x,y,z) \in [0,1] \times D_j.$$

We differentiate it with respect to $t$ and we get

$$0 = \int_{x_0}^{x} \int_{x_0}^{x_1} \int_{y_0}^{y} \int_{y_0}^{y_1} \int_{z_0}^{z} \int_{z_0}^{z_1} u_t(t,\alpha,\beta,\gamma) d\gamma dz_1 d\beta dy_1 d\alpha dx_1$$

$$+\int_{x_0}^{x} \int_{y_0}^{y} \int_{y_0}^{y_1} \int_{z_0}^{z} \int_{z_0}^{z_1} u^2(t,\alpha,\beta,\gamma) d\gamma dz_1 d\beta dy_1 d\alpha$$

$$+\int_{x_0}^{x} \int_{x_0}^{x_1} \int_{y_0}^{y} \int_{z_0}^{z} \int_{z_0}^{z_1} u(t,\alpha,\beta,\gamma) v(t,\alpha,\beta,\gamma) d\gamma dz_1 d\beta d\alpha dx_1$$

$$+\int_{x_0}^{x} \int_{x_0}^{x_1} \int_{y_0}^{y} \int_{y_0}^{y_1} \int_{z_0}^{z} u(t,\alpha,\beta,\gamma) w(t,\alpha,\beta,\gamma) d\gamma d\beta dy_1 d\alpha dx_1$$

$$+\int_{x_0}^{x} \int_{y_0}^{y} \int_{y_0}^{y_1} \int_{z_0}^{z} \int_{z_0}^{z_1} p(t,\alpha,\beta,\gamma) d\gamma dz_1 d\beta dy_1 d\alpha$$

$$-\int_{y_0}^{y} \int_{y_0}^{y_1} \int_{z_0}^{z} \int_{z_0}^{z_1} u(t,x,\beta,\gamma) d\gamma dz_1 d\beta dy_1$$

$$-\int_{x_0}^{x} \int_{x_0}^{x_1} \int_{z_0}^{z} \int_{z_0}^{z_1} u(t,\alpha,y,\gamma) d\gamma dz_1 d\alpha dx_1$$

$$-\int_{x_0}^{x} \int_{x_0}^{x_1} \int_{y_0}^{y} \int_{y_0}^{y_1} u(t,\alpha,\beta,z) d\beta dy_1 d\alpha dx_1, \qquad (t,x,y,z) \in [0,1] \times D_j.$$

We differentiate the last equality with respect to $x$ and we obtain

$$0 = \int_{x_0}^{x} \int_{y_0}^{y} \int_{y_0}^{y_1} \int_{z_0}^{z} \int_{z_0}^{z_1} u_t(t,\alpha,\beta,\gamma) d\gamma dz_1 d\beta dy_1 d\alpha$$

$$+\int_{y_0}^{y} \int_{y_0}^{y_1} \int_{z_0}^{z} \int_{z_0}^{z_1} u^2(t,x,\beta,\gamma) d\gamma dz_1 d\beta dy_1$$

$$+\int_{x_0}^{x}\int_{y_0}^{y}\int_{z_0}^{z}\int_{z_0}^{z_1} u(t,\alpha,\beta,\gamma)v(t,\alpha,\beta,\gamma)d\gamma dz_1 d\beta d\alpha$$

$$+\int_{x_0}^{x}\int_{y_0}^{y}\int_{y_0}^{y_1}\int_{z_0}^{z} u(t,\alpha,\beta,\gamma)w(t,\alpha,\beta,\gamma)d\gamma d\beta dy_1 d\alpha$$

$$+\int_{y_0}^{y}\int_{y_0}^{y_1}\int_{z_0}^{z}\int_{z_0}^{z_1} p(t,x,\beta,\gamma)d\gamma dz_1 d\beta dy_1$$

$$-\int_{y_0}^{y}\int_{y_0}^{y_1}\int_{z_0}^{z}\int_{z_0}^{z_1} u_x(t,x,\beta,\gamma)d\gamma dz_1 d\beta dy_1$$

$$-\int_{x_0}^{x}\int_{z_0}^{z}\int_{z_0}^{z_1} u(t,\alpha,y,\gamma)d\gamma dz_1 d\alpha$$

$$-\int_{x_0}^{x}\int_{y_0}^{y}\int_{y_0}^{y_1} u(t,\alpha,\beta,z)d\beta dy_1 d\alpha, (t,x,y,z)\in[0,1]\times D_j.$$

Again we differentiate with respect to $x$ and we get

$$0=\int_{y_0}^{y}\int_{y_0}^{y_1}\int_{z_0}^{z}\int_{z_0}^{z_1} u_t(t,x,\beta,\gamma)d\gamma dz_1 d\beta dy_1$$

$$+\int_{y_0}^{y}\int_{y_0}^{y_1}\int_{z_0}^{z}\int_{z_0}^{z_1} (u^2(t,x,\beta,\gamma))_x d\gamma dz_1 d\beta dy_1$$

$$+\int_{y_0}^{y}\int_{z_0}^{z}\int_{z_0}^{z_1} u(t,x,\beta,\gamma)v(t,x,\beta,\gamma)d\gamma dz_1 d\beta$$

$$+\int_{y_0}^{y}\int_{y_0}^{y_1}\int_{z_0}^{z} u(t,x,\beta,\gamma)w(t,x,\beta,\gamma)d\gamma d\beta dy_1$$

$$+\int_{y_0}^{y}\int_{y_0}^{y_1}\int_{z_0}^{z}\int_{z_0}^{z_1} p_x(t,x,\beta,\gamma)d\gamma dz_1 d\beta dy_1$$

$$-\int_{y_0}^{y}\int_{y_0}^{y_1}\int_{z_0}^{z}\int_{z_0}^{z_1} u_{xx}(t,x,\beta,\gamma)d\gamma dz_1 d\beta dy_1$$

$$-\int_{z_0}^{z}\int_{z_0}^{z_1} u(t,x,y,\gamma)d\gamma dz_1-\int_{y_0}^{y}\int_{y_0}^{y_1} u(t,x,\beta,z)d\beta dy_1, (t,x,y,z)\in[0,1]\times D_j.$$

Now we differentiate twice the last equation with respect to $y$ and we find

$$0=\int_{z_0}^{z}\int_{z_0}^{z_1} u_t(t,x,y,\gamma)d\gamma dz_1+\int_{z_0}^{z}\int_{z_0}^{z_1} (u^2(t,x,y,\gamma))_x d\gamma dz_1$$

$$+\int_{z_0}^{z}\int_{z_0}^{z_1} (u(t,x,y,\gamma)v(t,x,y,\gamma))_y d\gamma dz_1$$

$$+\int_{z_0}^{z} u(t,x,y,\gamma)w(t,x,y,\gamma)d\gamma+\int_{z_0}^{z}\int_{z_0}^{z_1} p_x(t,x,y,\gamma)d\gamma dz_1$$

$$-\int_{z_0}^{z}\int_{z_0}^{z_1} u_{yy}(t,x,y,\gamma)d\gamma dz_1-u(t,x,y,z)$$

$$-\int_{z_0}^{z}\int_{z_0}^{z_1} u_{xx}(t,x,y,\gamma)d\gamma dz_1, (t,x,y,z)\in[0,1]\times D_j.$$

We differentiate twice with respect to $z$ the last equation and we get

$$0=u_t(t,x,y,z)+(u^2(t,x,y,z))_x+(u(t,x,y,z)v(t,x,y,z))_y$$

$$+(u(t,x,y,z)w(t,x,y,z))_z + p_x(t,x,y,z)$$

$$-u_{xx}(t,x,y,z) - u_{yy}(t,x,y,z) - u_{zz}(t,x,y,z), (t,x,y,z) \in [0,1] \times D_j,$$

i.e., we get the first equation of the system (3.2).

As in above, after we differentiate with respect to $t$ and differentiate twice with respect to $x$, $y$ and $z$ the equations

$$I_2^{1j}(u,v,w,p) = 0, I_3^{1j}(u,v,w,p) = 0, I_4^{1j}(u,v,w,p) = 0$$

we get the second, third and fourth equation of (3.2), respectively.

After we put $t = 0$ in $I_1^{1j} = 0$, we obtain

$$0 = \int_{x_0}^{x} \int_{x_0}^{x_1} \int_{y_0}^{y} \int_{y_0}^{y_1} \int_{z_0}^{z} \int_{z_0}^{z_1} (u(0,\alpha,\beta,\gamma) - u_0(\alpha,\beta,\gamma)) d\gamma dz_1 d\beta dy_1 d\alpha dx_1,$$

which we differentiate twice in $x$, $y$ and $z$ and we find

$$u(0,x,y,z) = u_0(x,y,z), \qquad (x,y,z) \in D_j.$$

After we put $t = 0$ in $I_2^{1j} = 0$ and differentiate twice in $x$, $y$ and $z$ the equation $I_2^{1j} = 0$, we obtain

$$v(0,x,y,z) = v_0(x,y,z), \qquad (x,y,z) \in D_j.$$

After we put $t = 0$ in $I_3^{1j} = 0$ and differentiate twice in $x$, $y$ and $z$ the equation $I_3^{1j} = 0$, we obtain

$$w(0,x,y,z) = w_0(x,y,z), \qquad (x,y,z) \in D_j.$$

This completes the proof.

The proof of the existence result is based on Theorem 2.3.

Let $\tilde{\tilde{\tilde{X}}}^1$ be the set of all equicontinuous families of functions of the space

$$E^1 = \{g \in C^1([0,1], C_0^2(D_j)) : supp_{(x,y,z)} g \subset D_{jj}\}$$

with respect to the norm

$$||f|| = \max\{ \sup_{t\in[0,1],(x,y,z)\in D_j} |f(t,x,y,z)|,$$

$$\sup_{t\in[0,1],(x,y,z)\in D_j} |f_t(t,x,y,z)|,$$

$$\sup_{t\in[0,1],(x,y,z)\in D_j} |f_x(t,x,y,z)|, \qquad \sup_{t\in[0,1],(x,y,z)\in D_j} |f_{xx}(t,x,y,z)|,$$

$$\sup_{t\in[0,1],(x,y,z)\in D_j} |f_y(t,x,y,z)|, \qquad \sup_{t\in[0,1],(x,y,z)\in D_j} |f_{yy}(t,x,y,z)|,$$

$$\sup_{t\in[0,1],(x,y,z)\in D_j} |f_z(t,x,y,z)|, \qquad \sup_{t\in[0,1],(x,y,z)\in D_j} |f_{zz}(t,x,y,z)|\},$$

and $\tilde{\tilde{X}}^1 = \tilde{\tilde{\tilde{X}}}^1 \cup \{u_0, v_0, w_0\}$, $\tilde{X}^1 = \overline{\tilde{\tilde{X}}^1}$, i.e., $\tilde{X}^1$ is the completion of $\tilde{\tilde{X}}^1$, and

$$X^1 = \left\{ f \in \tilde{X}^1 : ||f|| \le M_{1j} := \frac{1}{2^{j2}\sqrt{\mu(D_j)}} \right\}.$$

Here $\mu(D_j)$ is the measure of the set $D_j$. Let also,

$$N_{1j} := \max\{\sup_{D_j}|u_0|, \qquad \sup_{D_j}|v_0|, \qquad \sup_{D_j}|w_0|\}$$

We take $\varepsilon > 0$ so that

$$\varepsilon\big(3M_{1j}^2 + 6M_{1j} + N_{1j}\big)\big(\mu(D_j)\big)^2 \leq M_{1j}.$$

We set

$$Y^1 = \big\{f \in \tilde{X}^1 : \| f \| \leq (1+\varepsilon)M_{1j}\big\}$$

By the construction of $X^1$ and $Y^1$, we have that $X^1$ is a compact subset of $Y^1$ and $Y^1$ is a compact subset of $C^1((0,1], C_0^2(D_j))$.

For $u, v, w, p \in Y^1$ we define the operators

$$S_1^{1j}(u,v,w,p) = -\varepsilon u + \varepsilon I_1^{1j}, \qquad T_1^{1j}(u,v,w,p) = (1+\varepsilon)u,$$

$$S_2^{1j}(u,v,w,p) = -\varepsilon v + \varepsilon I_2^{1j}, \qquad T_2^{1j}(u,v,w,p) = (1+\varepsilon)v,$$

$$S_3^{1j}(u,v,w,p) = -\varepsilon w + \varepsilon I_3^{1j}, \qquad T_3^{1j}(u,v,w,p) = (1+\varepsilon)w,$$

$$S_4^{1j}(u,v,w,p) = -\varepsilon p + \varepsilon I_4^{1j}, \qquad T_4^{1j}(u,v,w,p) = (1+\varepsilon)p,$$

$$S^{1j} = (S_1^{1j}, S_2^{1j}, S_3^{1j}, S_4^{1j}), \qquad T^{1j} = (T_1^{1j}, T_2^{1j}, T_3^{1j}, T_4^{1j}).$$

For $u, v, w, p \in X^1$ we have that

$$||I_1^{1j}(u,v,w,p)|| \leq (||u|| + ||u||^2 + ||u|| \cdot ||v|| + ||u|| \cdot ||w||$$

$$+||p|| + ||u|| + ||u|| + ||u|| + N_{1j})\big(\mu(D_j)\big)^2$$

$$\leq (3M_{1j}^2 + 5M_{1j} + N_{1j})\big(\mu(D_j)\big)^2.$$

Therefore, using our choice of $\varepsilon$

$$||S_1^{1j}|| \leq \varepsilon||u|| + \varepsilon(3M_{1j}^2 + 5M_{1j} + N_{1j})\big(\mu(D_j)\big)^2$$

$$\leq \varepsilon M_{1j} + \varepsilon(3M_{1j}^2 + 5M_{1j} + N_{1j})\big(\mu(D_j)\big)^2$$

$$\leq (1+\varepsilon)M_{1j}.$$

As in above we have

$$||S_2^{1j}|| \leq \varepsilon||v|| + \varepsilon(3M_{1j}^2 + 5M_{1j} + N_{1j})\big(\mu(D_j)\big)^2$$

$$\leq \varepsilon M_{1j} + \varepsilon(3M_{1j}^2 + 5M_{1j} + N_{1j})\big(\mu(D_j)\big)^2$$

$$\leq (1+\varepsilon)M_{1j},$$

$$||S_3^{1j}|| \leq \varepsilon||w|| + \varepsilon(3M_{1j}^2 + 5\mathrm{M}_{1j} + N_{1j})\big(\mu(D_j)\big)^2$$

$$\leq \varepsilon M_{1j} + \varepsilon(3M_{1j}^2 + 5M_{1j} + N_{1j})\big(\mu(D_j)\big)^2$$

$$\leq (1+\varepsilon)M_{1j},$$

$$||S_4^{1j}|| \leq \varepsilon||p|| + 3\varepsilon M_{1j}\big(\mu(D_j)\big)^2$$

$$\leq \varepsilon M_{1j} + M_{1j}$$

$$= (1+\varepsilon)M_{1j}.$$

Therefore, for $(u, v, w, p) \in X^1$ we have that

$$S_i^{1j}(u, v, w, p) \in Y^1, \qquad i = 1,2,3,4.$$

Then

$$S^{1j}: X^1 \times X^1 \times X^1 \times X^1 \to Y^1 \times Y^1 \times Y^1 \times Y^1$$

and it is continuous.

The operator

$$T^{1j}: X^1 \times X^1 \times X^1 \times X^1 \to Y^1 \times Y^1 \times Y^1 \times Y^1$$

is an expansive operator with constant $1+\varepsilon > 1$ and if $(u, v, w, p) \in Y^1 \times Y^1 \times Y^1 \times Y^1$, then

$$\left(\frac{1}{1+\varepsilon}u, \frac{1}{1+\varepsilon}v, \frac{1}{1+\varepsilon}w, \frac{1}{1+\varepsilon}p\right) \in X^1 \times X^1 \times X^1 \times X^1,$$

and

$$(T_1^{1j}\left(\frac{1}{1+\varepsilon}u, \frac{1}{1+\varepsilon}v, \frac{1}{1+\varepsilon}w, \frac{1}{1+\varepsilon}p\right), T_2^{1j}\left(\frac{1}{1+\varepsilon}u, \frac{1}{1+\varepsilon}v, \frac{1}{1+\varepsilon}w, \frac{1}{1+\varepsilon}p\right),$$

$$T_3^{1j}\left(\frac{1}{1+\varepsilon}u, \frac{1}{1+\varepsilon}v, \frac{1}{1+\varepsilon}w, \frac{1}{1+\varepsilon}p\right), T_4^{1j}\left(\frac{1}{1+\varepsilon}u, \frac{1}{1+\varepsilon}v, \frac{1}{1+\varepsilon}w, \frac{1}{1+\varepsilon}p\right))$$

$$= (u, v, w, p).$$

Consequently $T^{1j}: X^1 \times X^1 \times X^1 \times X^1 \to Y^1 \times Y^1 \times Y^1 \times Y^1$ is onto.

From here and from Theorem 2.3, it follows that the operator $T^{1j} + S^{1j}$ has a fixed point $(u_1, v_1, w_1, p_1)$ in $X^1 \times X^1 \times X^1 \times X^1$. For it we have

$$T_1^{1j}(u_1, v_1, w_1, p_1) + S_1^{1j}(u_1, v_1, w_1, p_1) = u_1$$

$$T_2^{1j}(u_1, v_1, w_1, p_1) + S_2^{1j}(u_1, v_1, w_1, p_1) = \mathrm{v}_1$$

$$T_3^{1j}(u_1, v_1, w_1, p_1) + S_3^{1j}(u_1, v_1, w_1, p_1) = w_1$$

$$T_4^{1j}(u_1, v_1, w_1, p_1) + S_4^{1j}(u_1, v_1, w_1, p_1) = p_1$$

or

$$(1+\varepsilon)u_1 - \varepsilon u_1 + I_1^{1j}(u_1, v_1, w_1, p_1) = u_1$$

$$(1+\varepsilon)v_1 - \varepsilon v_1 + I_2^{1j}(u_1, v_1, w_1, p_1) = v_1$$

$$(1+\varepsilon)w_1 - \varepsilon w_1 + I_3^{1j}(u_1, v_1, w_1, p_1) = w_1$$

$$(1+\varepsilon)p_1 - \varepsilon p_1 + I_4^{1j}(u_1, v_1, w_1, p_1) = p_1,$$

whereupon

$$I_1^{1j}(u_1,v_1,w_1,p_1)=0, \qquad I_2^{1j}(u_1,v_1,w_1,p_1)=0,$$

$$I_3^{1j}(u_1,v_1,w_1,p_1)=0, \qquad I_4^{1j}(u_1,v_1,w_1,p_1)=0.$$

Hence and Lemma 3.2 we obtain that $(u_1,v_1,w_1,p_1)$ is a solution of the system (3.2) for which $u_1,\ v_1,\ w_1,\ p_1\in C^1([0,1],C_0^2(D_j))$.

Figure 1: A sketch of the division of $\mathbb{R}^3$ to disjointed subsets D1,D2,D3 etc, that illustrates how the proof of existence steps are done.

**Remark 3.3** *Since in the estimates of* $\parallel I_i^{1j}(u,v,w,p)\parallel$ *participate* $\mu(D_j)$*, we suggest* $D_j$ *to be bounded,* $j\in\mathbb{N}$*.*

If we assume that

$$u_1(t,x,y,z)=u_0(x,y,z),$$

$$v_1(t,x,y,z)=v_0(x,y,z),$$

$$w_1(t,x,y,z)=w_0(x,y,z), \qquad t\in[0,1], \qquad (x,y,z)\in D_j,$$

then, using $I_l^{ij}(u_1,v_1,w_1,p)(t,x,y,z)=0,\ t\in[0,1],\ (x,y,z)\in D_j,\ l=1,2,3,$ we get

$$\begin{aligned}
& t\int_{x_0}^{x}\int_{y_0}^{y}\int_{y_0}^{y_1}\int_{z_0}^{z}\int_{z_0}^{z_1}(u_0(\alpha,\beta,\gamma))^2 d\gamma dz_1 d\beta dy_1 d\alpha dx_1 \\
& +t\int_{x_0}^{x}\int_{x_0}^{x_1}\int_{y_0}^{y}\int_{z_0}^{z}\int_{z_0}^{z_1}u_0(\alpha,\beta,\gamma)v_0(\alpha,\beta,\gamma)d\gamma dz_1 d\beta d\alpha dx_1 \\
& +t\int_{x_0}^{x}\int_{x_0}^{x_1}\int_{y_0}^{y}\int_{y_0}^{y_1}\int_{z_0}^{z}u_0(\alpha,\beta,\gamma)w_0(\alpha,\beta,\gamma)d\gamma d\beta dy_1 d\alpha dx_1 \\
& +\int_{0}^{t}\int_{x_0}^{x}\int_{y_0}^{y}\int_{y_0}^{y_1}\int_{z_0}^{z}\int_{z_0}^{z_1}p_1(s,\alpha,\beta,\gamma)d\gamma dz_1 d\beta dy_1 d\alpha ds \qquad (3.4) \\
& -t\int_{y_0}^{y}\int_{y_0}^{y_1}\int_{z_0}^{z}\int_{z_0}^{z_1}u_0(x,\beta,\gamma)d\gamma dz_1 d\beta dy_1 \\
& -t\int_{x_0}^{x}\int_{x_0}^{x_1}\int_{y_0}^{y}\int_{y_0}^{y_1}u_0(\alpha,\beta,z)d\beta dy_1 d\alpha dx_1 \\
& -t\int_{x_0}^{x}\int_{x_0}^{x_1}\int_{z_0}^{z}\int_{z_0}^{z_1}u_0(\alpha,y,\gamma)d\gamma dz_1 d\alpha dx_1=0, \qquad t\in[0,1],(x,y,z)\in D_j,
\end{aligned}$$

$$
\begin{aligned}
& t\int_{x_0}^{x}\int_{y_0}^{y}\int_{y_0}^{y_1}\int_{z_0}^{z}\int_{z_0}^{z_1} u_0(\alpha,\beta,\gamma)v_0(\alpha,\beta,\gamma)d\gamma dz_1 d\beta dy_1 d\alpha dx_1 \\
& +t\int_{x_0}^{x}\int_{x_0}^{x_1}\int_{y_0}^{y}\int_{z_0}^{z}\int_{z_0}^{z_1} (v_0(\alpha,\beta,\gamma))^2 d\gamma dz_1 d\beta d\alpha dx_1 \\
& +t\int_{x_0}^{x}\int_{x_0}^{x_1}\int_{y_0}^{y}\int_{y_0}^{y_1}\int_{z_0}^{z} v_0(\alpha,\beta,\gamma)w_0(\alpha,\beta,\gamma)d\gamma d\beta dy_1 d\alpha dx_1 \\
& +\int_{0}^{t}\int_{x_0}^{x}\int_{x_0}^{x_1}\int_{y_0}^{y}\int_{z_0}^{z}\int_{z_0}^{z_1} p_1(s,\alpha,\beta,\gamma)d\gamma dz_1 d\beta d\alpha dx_1 ds \\
& -t\int_{y_0}^{y}\int_{y_0}^{y_1}\int_{z_0}^{z}\int_{z_0}^{z_1} v_0(x,\beta,\gamma)d\gamma dz_1 d\beta dy_1 \\
& -t\int_{x_0}^{x}\int_{x_0}^{x_1}\int_{y_0}^{y}\int_{y_0}^{y_1} v_0(\alpha,\beta,z)d\beta dy_1 d\alpha dx_1 \\
& -t\int_{x_0}^{x}\int_{x_0}^{x_1}\int_{z_0}^{z}\int_{z_0}^{z_1} v_0(\alpha,y,\gamma)d\gamma dz_1 d\alpha dx_1 = 0, \qquad t\in[0,1],(x,y,z)\in D_j,
\end{aligned}
\tag{3.5}
$$

$$
\begin{aligned}
& t\int_{x_0}^{x}\int_{y_0}^{y}\int_{y_0}^{y_1}\int_{z_0}^{z}\int_{z_0}^{z_1} u_0(\alpha,\beta,\gamma)w_0(\alpha,\beta,\gamma)d\gamma dz_1 d\beta dy_1 d\alpha dx_1 \\
& +t\int_{x_0}^{x}\int_{x_0}^{x_1}\int_{y_0}^{y}\int_{z_0}^{z}\int_{z_0}^{z_1} w_0(\alpha,\beta,\gamma)v_0(\alpha,\beta,\gamma)d\gamma dz_1 d\beta d\alpha dx_1 \\
& +t\int_{x_0}^{x}\int_{x_0}^{x_1}\int_{y_0}^{y}\int_{y_0}^{y_1}\int_{z_0}^{z} (w_0(\alpha,\beta,\gamma))^2 d\gamma d\beta dy_1 d\alpha dx_1 \\
& +\int_{0}^{t}\int_{x_0}^{x}\int_{x_0}^{x_1}\int_{y_0}^{y}\int_{y_0}^{y_1}\int_{z_0}^{z} p_1(s,\alpha,\beta,\gamma)d\gamma dd\beta dy_1 d\alpha dx_1 ds \\
& -t\int_{y_0}^{y}\int_{y_0}^{y_1}\int_{z_0}^{z}\int_{z_0}^{z_1} w_0(x,\beta,\gamma)d\gamma dz_1 d\beta dy_1 \\
& -t\int_{x_0}^{x}\int_{x_0}^{x_1}\int_{y_0}^{y}\int_{y_0}^{y_1} w_0(\alpha,\beta,z)d\beta dy_1 d\alpha dx_1 \\
& -t\int_{x_0}^{x}\int_{x_0}^{x_1}\int_{z_0}^{z}\int_{z_0}^{z_1} w_0(\alpha,y,\gamma)d\gamma dz_1 d\alpha dx_1 = 0, \qquad t\in[0,1],(x,y,z)\in D_j,
\end{aligned}
\tag{3.6}
$$

$$
\begin{aligned}
& t\int_{x_0}^{x}\int_{y_0}^{y}\int_{y_0}^{y_1}\int_{z_0}^{z}\int_{z_0}^{z_1} u_0(\alpha,\beta,\gamma)d\gamma dz_1 d\beta dy_1 d\alpha \\
& +t\int_{x_0}^{x}\int_{x_0}^{x_1}\int_{y_0}^{y}\int_{y_0}^{y_1}\int_{z_0}^{z} u_0(\alpha,\beta,\gamma)d\gamma d\beta dy_1 d\alpha dx_1 \\
& +t\int_{x_0}^{x}\int_{x_0}^{x_1}\int_{y_0}^{y}\int_{z_0}^{z}\int_{z_0}^{z_1} u_0(\alpha,\beta,\gamma)d\gamma dz_1 d\beta d\alpha dx_1 = 0, \qquad t\in[0,1],(x,y,z)\in D_j.
\end{aligned}
\tag{3.7}
$$

If the initial functions $u_0$, $v_0$, $w_0$ and $p_1$ do not satisfy (3.4), (3.5), (3.6), (3.7), then

$$(u_0,v_0,w_0)\not\cong(u_1,v_1,w_1) \qquad on \qquad [0,1]\times D_j.$$

Now we consider the problem

$$u_t + (u^2)_x + (uv)_y + (uw)_z + p_x - u_{xx} - u_{yy} - u_{zz} = 0$$

$$v_t + (uv)_x + (v^2)_y + (vw)_z + p_y - v_{xx} - v_{yy} - v_{zz} = 0$$

$$w_t + (uw)_x + (vw)_y + (w^2)_z + p_z - w_{xx} - w_{yy} - w_{zz} = 0 \tag{3.8}$$

$$u_x + v_y + w_z = 0 \qquad in \qquad (1,2] \times D_j,$$

$$u(1,x,y,z) = u_1(1,x,y,z), \qquad v(1,x,y,z) = v_1(1,x,y,z),$$

$$w(1,x,y,z) = w_1(1,x,y,z), \qquad (x,y,z) \in D_j.$$

We will prove that the problem (3.8) has a solution $(u,v,w,p)$ such that $u,v,w,p \in C^1([1,2], C_0^2(D_j))$.

For $(u,v,w,p) \in \left(C^1([1,2], C_0^2(D_j))\right)^4$ we define

$$I_1^{2j}(u,v,w,p) = \int_{x_0}^{x}\int_{x_0}^{x_1}\int_{y_0}^{y}\int_{y_0}^{y_1}\int_{z_0}^{z}\int_{z_0}^{z_1}(u(t,\alpha,\beta,\gamma)$$

$$-u_1(1,\alpha,\beta,\gamma))d\gamma dz_1 d\beta dy_1 d\alpha dx_1$$

$$+\int_1^t\int_{x_0}^{x}\int_{y_0}^{y}\int_{y_0}^{y_1}\int_{z_0}^{z}\int_{z_0}^{z_1}u^2(s,\alpha,\beta,\gamma)d\gamma dz_1 d\beta dy_1 d\alpha ds$$

$$+\int_1^t\int_{x_0}^{x}\int_{x_0}^{x_1}\int_{y_0}^{y}\int_{z_0}^{z}\int_{z_0}^{z_1}u(s,\alpha,\beta,\gamma)v(s,\alpha,\beta,\gamma)d\gamma dz_1 d\beta d\alpha dx_1 ds$$

$$+\int_1^t\int_{x_0}^{x}\int_{x_0}^{x_1}\int_{y_0}^{y}\int_{y_0}^{y_1}\int_{z_0}^{z}u(s,\alpha,\beta,\gamma)w(s,\alpha,\beta,\gamma)d\gamma d\beta dy_1 d\alpha dx_1 ds$$

$$+\int_1^t\int_{x_0}^{x}\int_{y_0}^{y}\int_{y_0}^{y_1}\int_{z_0}^{z}\int_{z_0}^{z_1}p(s,\alpha,\beta,\gamma)d\gamma dz_1 d\beta dy_1 d\alpha ds$$

$$-\int_1^t\int_{y_0}^{y}\int_{y_0}^{y_1}\int_{z_0}^{z}\int_{z_0}^{z_1}u(s,x,\beta,\gamma)d\gamma dz_1 d\beta dy_1 ds$$

$$-\int_1^t\int_{x_0}^{x}\int_{x_0}^{x_1}\int_{z_0}^{z}\int_{z_0}^{z_1}u(s,\alpha,y,\gamma)d\gamma dz_1 d\alpha dx_1 ds$$

$$-\int_1^t\int_{x_0}^{x}\int_{x_0}^{x_1}\int_{y_0}^{y}\int_{y_0}^{y_1}u(s,\alpha,\beta,z)d\beta dy_1 d\alpha dx_1 ds,$$

$$I_2^{2j}(u,v,w,p) = \int_{x_0}^{x}\int_{x_0}^{x_1}\int_{y_0}^{y}\int_{y_0}^{y_1}\int_{z_0}^{z}\int_{z_0}^{z_1}(v(t,\alpha,\beta,\gamma)$$

$$-v_1(1,\alpha,\beta,\gamma))d\gamma dz_1 d\beta dy_1 d\gamma dx_1$$

$$+\int_1^t\int_{x_0}^{x}\int_{y_0}^{y}\int_{y_0}^{y_1}\int_{z_0}^{z}\int_{z_0}^{z_1}u(s,\alpha,\beta,\gamma)v(s,\alpha,\beta,\gamma)d\gamma dz_1 d\beta dy_1 d\alpha ds$$

$$+\int_1^t \int_{x_0}^{x} \int_{x_0}^{x_1} \int_{y_0}^{y} \int_{z_0}^{z} \int_{z_0}^{z_1} v^2(s,\alpha,\beta,\gamma) d\gamma dz_1 d\beta d\alpha dx_1 ds$$

$$+\int_1^t \int_{x_0}^{x} \int_{x_0}^{x_1} \int_{y_0}^{y} \int_{y_0}^{y_1} \int_{z_0}^{z} v(s,\alpha,\beta,\gamma) w(s,\alpha,\beta,\gamma) d\gamma d\beta dy_1 d\alpha dx_1 ds$$

$$+\int_1^t \int_{x_0}^{x} \int_{x_0}^{x_1} \int_{y_0}^{y} \int_{z_0}^{z} \int_{z_0}^{z_1} p(s,\alpha,\beta,\gamma) d\gamma dz_1 d\beta d\alpha dx_1 ds$$

$$-\int_1^t \int_{y_0}^{y} \int_{y_0}^{y_1} \int_{z_0}^{z} \int_{z_0}^{z_1} v(s,x,\beta,\gamma) d\gamma dz_1 d\beta dy_1 ds$$

$$-\int_1^t \int_{x_0}^{x} \int_{x_0}^{x_1} \int_{z_0}^{z} \int_{z_0}^{z_1} v(s,\alpha,y,\gamma) d\gamma dz_1 d\alpha dx_1 ds$$

$$-\int_1^t \int_{x_0}^{x} \int_{x_0}^{x_1} \int_{y_0}^{y} \int_{y_0}^{y_1} v(s,\alpha,\beta,z) d\beta dy_1 d\alpha dx_1 ds,$$

$$I_3^{2j}(u,v,w,p) = \int_{x_0}^{x} \int_{x_0}^{x_1} \int_{y_0}^{y} \int_{y_0}^{y_1} \int_{z_0}^{z} \int_{z_0}^{z_1} (w(t,\alpha,\beta,\gamma)$$

$$-w_1(1,\alpha,\beta,\gamma)) d\gamma dz_1 d\beta dy_1 d\gamma dx_1$$

$$+\int_1^t \int_{x_0}^{x} \int_{y_0}^{y} \int_{y_0}^{y_1} \int_{z_0}^{z} \int_{z_0}^{z_1} u(s,\alpha,\beta,\gamma) w(s,\alpha,\beta,\gamma) d\gamma dz_1 d\beta dy_1 d\alpha ds$$

$$+\int_1^t \int_{x_0}^{x} \int_{x_0}^{x_1} \int_{y_0}^{y} \int_{z_0}^{z} \int_{z_0}^{z_1} v(s,\alpha,\beta,\gamma) w(s,\alpha,\beta,\gamma) d\gamma dz_1 d\beta d\alpha dx_1 ds$$

$$+\int_1^t \int_{x_0}^{x} \int_{x_0}^{x_1} \int_{y_0}^{y} \int_{y_0}^{y_1} \int_{z_0}^{z} w^2(s,\alpha,\beta,\gamma) d\gamma d\beta dy_1 d\alpha dx_1 ds$$

$$+\int_1^t \int_{x_0}^{x} \int_{x_0}^{x_1} \int_{y_0}^{y} \int_{y_0}^{y_1} \int_{z_0}^{z} p(s,\alpha,\beta,\gamma) d\gamma d\beta dy_1 d\alpha dx_1 ds$$

$$-\int_1^t \int_{y_0}^{y} \int_{y_0}^{y_1} \int_{z_0}^{z} \int_{z_0}^{z_1} w(s,x,\beta,\gamma) d\gamma dz_1 d\beta dy_1 ds$$

$$-\int_1^t \int_{x_0}^{x} \int_{x_0}^{x_1} \int_{z_0}^{z} \int_{z_0}^{z_1} w(s,\alpha,y,\gamma) d\gamma dz_1 d\alpha dx_1 ds$$

$$-\int_1^t \int_{x_0}^{x} \int_{x_0}^{x_1} \int_{y_0}^{y} \int_{y_0}^{y_1} w(s,\alpha,\beta,z) d\beta dy_1 d\alpha dx_1 ds,$$

$$I_4^{2j}(u,v,w,p) = \int_1^t \int_{x_0}^{x} \int_{y_0}^{y} \int_{y_0}^{y_1} \int_{z_0}^{z} \int_{z_0}^{z_1} u(s,\alpha,\beta,\gamma) d\gamma dz_1 d\beta dy_1 d\alpha ds$$

$$+\int_1^t \int_{x_0}^{x} \int_{x_0}^{x_1} \int_{y_0}^{y} \int_{z_0}^{z} \int_{z_0}^{z_1} v(s,\alpha,\beta,\gamma) d\gamma dz_1 d\beta d\alpha dx_1 ds$$

$$+\int_1^t \int_{x_0}^{x} \int_{x_0}^{x_1} \int_{y_0}^{y} \int_{y_0}^{y_1} \int_{z_0}^{z} w(s,\alpha,\beta,\gamma) d\gamma d\beta dy_1 d\alpha dx_1 ds.$$

We note that after we differentiate with respect to $t$ and twice with respect to $x$, $y$ and $z$ the system

$$I_1^{2j}(u,v,w,p)=0, \qquad I_2^{2j}(u,v,w,p)=0,$$
$$I_3^{2j}(u,v,w,p)=0, \qquad I_4^{2j}(u,v,w,p)=0, \tag{3.9}$$

we get the system (3.8). After we put $t=1$ in $I_1^{2j}=0$ and differentiate twice in $x$, $y$ and $z$ the equation $I_1^{2j}=0$ we obtain

$$u(1,x,y,z)=u_1(1,x,y,z), \qquad (x,y,z)\in D_j.$$

After we put $t=1$ in $I_2^{2j}=0$ and differentiate twice in $x$, $y$ and $z$ the equation $I_2^{2j}=0$ we obtain

$$v(1,x,y,z)=v_1(1,x,y,z), \qquad (x,y,z)\in D_j.$$

After we put $t=1$ in $I_3^{2j}=0$ and differentiate twice in $x$, $y$ and $z$ the equation $I_3^{2j}=0$ we obtain

$$w(1,x,y,z)=w_1(1,x,y,z), \qquad (x,y,z)\in D_j.$$

Consequently every solution $(u,v,w,p)\in\left(C^1([1,2],C_0^2(D_j))\right)^4$ of (3.9) is a solution of the problem (3.8).

Let $\tilde{\tilde{\tilde{X}}}^2$ be the set of all equicontinuous families of the space

$$\{g\in C^1([1,2],C_0^2(D_j)), supp_{x,y,z}g\subset D_{jj}\}$$

with respect to the norm

$$||f||=\max\{\sup_{t\in[1,2],(x,y,z)\in D_j}|f(t,x,y,z)|,$$

$$\sup_{t\in[1,2],(x,y,z)\in D_j}|f_t(t,x,y,z)|,$$

$$\sup_{t\in[1,2],(x,y,z)\in D_j}|f_x(t,x,y,z)|, \qquad \sup_{t\in[1,2],(x,y,z)\in D_j}|f_{xx}(t,x,y,z)|,$$

$$\sup_{t\in[1,2],(x,y,z)\in D_j}|f_y(t,x,y,z)|, \qquad \sup_{t\in[1,2],(x,y,z)\in D_j}|f_{yy}(t,x,y,z)|,$$

$$\sup_{t\in[1,2],(x,y,z)\in D_j}|f_z(t,x,y,z)|, \qquad \sup_{t\in[1,2],(x,y,z)\in D_j}|f_{zz}(t,x,y,z)|\},$$

$$f\in\tilde{X}^2,$$

and

$$\tilde{\tilde{X}}^2=\tilde{\tilde{\tilde{X}}}^2\cup\{u_1(1,x,y,z),v_1(1,x,y,z),w_1(1,x,y,z)\},$$

$\tilde{X}^2=\overline{\tilde{\tilde{X}}^2}$, i.e., $\tilde{X}^2$ is the completition of $\tilde{\tilde{X}}^2$, and

$$X^2=\left\{f\in\tilde{X}^2:\|f\|\leq M_j=\frac{1}{2^{j^2}\sqrt{\mu(D_j)}}\right\},$$

and

$$Y^2=\left\{f\in\tilde{X}^2:\|f\|\leq(1+\varepsilon)M_{1j}\right\}.$$

Note that $X^2$ is a compact subset of $Y^2$.

For $u,v,w,p\in Y^2$ we define the operators

$$S_1^{2j}(u,v,w,p) = -\varepsilon u + \varepsilon I_1^{2j}, \qquad T_1^{2j}(u,v,w,p) = (1+\varepsilon)u,$$

$$S_2^{2j}(u,v,w,p) = -\varepsilon v + \varepsilon I_2^{2j}, \qquad T_2^{2j}(u,v,w,p) = (1+\varepsilon)v,$$

$$S_3^{2j}(u,v,w,p) = -\varepsilon w + \varepsilon I_3^{2j}, \qquad T_3^{2j}(u,v,w,p) = (1+\varepsilon)w,$$

$$S_4^{2j}(u,v,w,p) = -\varepsilon p + \varepsilon I_4^{2}, \qquad T_4^{2j}(u,v,w,p) = (1+\varepsilon)p,$$

$$S^{2j} = (S_1^{2j}, S_2^{2j}, S_3^{2j}, S_4^{2j}), \qquad T^{2j} = (T_1^{2j}, T_2^{2j}, T_3^{2j}, T_4^{2j}).$$

For $u,v,w,p \in X^2$ we have that $S_i^{2j}(u,v,w,p) \in Y^2$, $i = 1,2,3,4$, i.e.,

$$S^{2j}: X^2 \times X^2 \times X^2 \times X^2 \to Y^2 \times Y^2 \times Y^2 \times Y^2$$

and it is continuous.

The operator

$$T^{2j}: X^2 \times X^2 \times X^2 \times X^2 \to Y^2 \times Y^2 \times Y^2 \times Y^2$$

is an expansive operator with constant $1+\varepsilon > 1$ and if $(u,v,w,p) \in Y^2 \times Y^2 \times Y^2 \times Y^2$, then

$$\left(\frac{1}{1+\varepsilon}u, \frac{1}{1+\varepsilon}v, \frac{1}{1+\varepsilon}w, \frac{1}{1+\varepsilon}p\right) \in X^2 \times X^2 \times X^2 \times X^2,$$

and

$$(T_1^{2j}\left(\frac{1}{1+\varepsilon}u, \frac{1}{1+\varepsilon}v, \frac{1}{1+\varepsilon}w, \frac{1}{1+\varepsilon}p\right), T_2^{2j}\left(\frac{1}{1+\varepsilon}u, \frac{1}{1+\varepsilon}v, \frac{1}{1+\varepsilon}w, \frac{1}{1+\varepsilon}p\right),$$

$$T_3^{2j}\left(\frac{1}{1+\varepsilon}u, \frac{1}{1+\varepsilon}v, \frac{1}{1+\varepsilon}w, \frac{1}{1+\varepsilon}p\right), T_4^{2j}\left(\frac{1}{1+\varepsilon}u, \frac{1}{1+\varepsilon}v, \frac{1}{1+\varepsilon}w, \frac{1}{1+\varepsilon}p\right))$$

$$= (u,v,w,p).$$

Consequently $T^{2j}: X^2 \times X^2 \times X^2 \times X^2 \to Y^2 \times Y^2 \times Y^2 \times Y^2$ is onto.

From here and from Theorem 2.3, it follows that the operator $T^{2j} + S^{2j}$ has a fixed point $(u_2, v_2, w_2, p_2)$ in $X^2 \times X^2 \times X^2 \times X^2$. For it we have

$$T_1^{2j}(u_2,v_2,w_2,p_2) + S_1^{2j}(u_2,v_2,w_2,p_2) = u_2$$

$$T_2^{2j}(u_2,v_2,w_2,p_2) + S_2^{2j}(u_2,v_2,w_2,p_2) = v_2$$

$$T_3^{2j}(u_2,v_2,w_2,p_2) + S_3^{2j}(u_2,v_2,w_2,p_2) = w_2$$

$$T_4^{2j}(u_2,v_2,w_2,p_2) + S_4^{2j}(u_2,v_2,w_2,p_2) = p_2$$

or

$$(1+\varepsilon)u_2 - \varepsilon u_2 + I_1^{2j}(u_2, v_2, w_2, p_2) = u_2$$

$$(1+\varepsilon)v_2 - \varepsilon v_2 + I_2^{2j}(u_2, v_2, w_2, p_2) = v_2$$

$$(1+\varepsilon)w_2 - \varepsilon w_2 + I_3^{2j}(u_2, v_2, w_2, p_2) = w_2$$

$$(1+\varepsilon)p_2 - \varepsilon p_2 + I_4^{2j}(u_2, v_2, w_2, p_2) = p_2,$$

whereupon

$$I_1^{2j}(u_2, v_2, w_2, p_2) = 0, \qquad I_2^{2j}(u_2, v_2, w_2, p_2) = 0,$$

$$I_3^{2j}(u_2, v_2, w_2, p_2) = 0, \qquad I_4^{2j}(u_2, v_2, w_2, p_2) = 0.$$

Therefore $(u_2, v_2, w_2, p_2)$ is a solution of the system (3.8) for which $u_2$, $v_2$, $w_2$, $p_2 \in C^1([1,2], C_0^2(D_j))$.

We note that

$$u_1(1,x,y,z) = u_2(1,x,y,z),$$

$$v_1(1,x,y,z) = v_2(1,x,y,z),$$

$$w_1(1,x,y,z) = w_2(1,x,y,z),$$

$$p_1(1,x,y,z) = p_2(1,x,y,z), \qquad (x,y,z) \in D_j,$$

whereupon

$$u_{1x}(1,x,y,z) = u_{2x}(1,x,y,z),$$

$$v_{1x}(1,x,y,z) = v_{2x}(1,x,y,z),$$

$$w_{1x}(1,x,y,z) = w_{2x}(1,x,y,z),$$

$$p_{1x}(1,x,y,z) = p_{2x}(1,x,y,z),$$

$$u_{1y}(1,x,y,z) = u_{2y}(1,x,y,z),$$

$$v_{1y}(1,x,y,z) = v_{2y}(1,x,y,z),$$

$$w_{1y}(1,x,y,z) = w_{2y}(1,x,y,z),$$

$$p_{1y}(1,x,y,z) = p_{2y}(1,x,y,z),$$

$$u_{1z}(1,x,y,z) = u_{2z}(1,x,y,z),$$

$$v_{1z}(1,x,y,z) = v_{2z}(1,x,y,z),$$

$$w_{1z}(1,x,y,z) = w_{2z}(1,x,y,z),$$

$$p_{1z}(1,x,y,z) = p_{2z}(1,x,y,z),$$

$$u_{1xx}(1,x,y,z) = u_{2xx}(1,x,y,z),$$

$$v_{1xx}(1,x,y,z) = v_{2xx}(1,x,y,z),$$

$$w_{1xx}(1,x,y,z) = w_{2xx}(1,x,y,z),$$

$$p_{1xx}(1,x,y,z) = p_{2xx}(1,x,y,z),$$

$$u_{1yy}(1,x,y,z) = u_{2yy}(1,x,y,z),$$

$$v_{1yy}(1,x,y,z) = v_{2yy}(1,x,y,z),$$

$$w_{1yy}(1,x,y,z) = w_{2yy}(1,x,y,z),$$

$$p_{1yy}(1,x,y,z) = p_{2yy}(1,x,y,z),$$

$$u_{1zz}(1,x,y,z) = u_{2zz}(1,x,y,z),$$

$$v_{1zz}(1,x,y,z) = v_{2zz}(1,x,y,z),$$

$$w_{1zz}(1,x,y,z) = w_{2zz}(1,x,y,z),$$

$$p_{1zz}(1,x,y,z) = p_{2zz}(1,x,y,z), \qquad (x,y,z) \in D_j.$$

Hence and (3.2), (3.8), we get

$$u_{1t}(1,x,y,z) = u_{2t}(1,x,y,z),$$

$$v_{1t}(1,x,y,z) = v_{2t}(1,x,y,z),$$

$$w_{1t}(1,x,y,z) = w_{2t}(1,x,y,z),$$

$$p_{1t}(1,x,y,z) = p_{2t}(1,x,y,z), \qquad (x,y,z) \in D_j.$$

Consequently

$$(u(t,x,y,z), v(t,x,y,z), w(t,x,y,z), p(t,x,y,z))$$

$$= \begin{cases} (u_1(t,x,y,z), v_1(t,x,y,z), w_1(t,x,y,z), p_1(t,x,y,z)) \in (C^1([0,1], C_0^2(D_j)))^4 \\ \\ (u_2(t,x,y,z), v_2(t,x,y,z), w_2(t,x,y,z), p_2(t,x,y,z)) \in (C^1([1,2], C_0^2(D_j)))^4 \end{cases}$$

belongs to $(C^1([0,2], C_0^2(D_j)))^4$ and it is a solution to the problem

$$u_t + (u^2)_x + (uv)_y + (uw)_z + p_x - u_{xx} - u_{yy} - u_{zz} = 0$$

$$v_t + (uv)_x + (v^2)_y + (vw)_z + p_y - v_{xx} - v_{yy} - v_{zz} = 0$$

$$w_t + (uw)_x + (vw)_y + (w^2)_z + p_z - w_{xx} - w_{yy} - w_{zz} = 0$$

$$u_x + v_y + w_z = 0 \qquad in \qquad (0,2] \times D_j,$$

$$u(0,x,y,z) = u_0(x,y,z), \qquad v(0,x,y,z) = v_0(x,y,z),$$

$$w(0,x,y,z) = w_0(x,y,z), \qquad (x,y,z) \in D_j.$$

Then we consider the problem

$$u_t + (u^2)_x + (uv)_y + (uw)_z + p_x - u_{xx} - u_{yy} - u_{zz} = 0$$

$$v_t + (uv)_x + (v^2)_y + (vw)_z + p_y - v_{xx} - v_{yy} - v_{zz} = 0$$

$$w_t + (uw)_x + (vw)_y + (w^2)_z + p_z - w_{xx} - w_{yy} - w_{zz} = 0$$

$$u_x + v_y + w_z = 0 \qquad in \qquad (2,3] \times D_j,$$

$$u(2,x,y,z) = u_2(2,x,y,z), \qquad v(2,x,y,z) = v_2(2,x,y,z),$$

$$w(2,x,y,z) = w_2(2,x,y,z), \qquad (x,y,z) \in D_j$$

and as above we construct a solution

$$(u_3, v_3, w_3, p_3) \in \left(C^1([2,3], C_0^2(D_j))\right)^4$$

and so on. Consequently

$$(u^j(t,x,y,z),v^j(t,x,y,z),w^j(t,x,y,z),p^j(t,x,y,z))$$

$$= \begin{cases} (u_1(t,x,y,z),v_1(t,x,y,z),w_1(t,x,y,z),p_1(t,x,y,z)) \in (C^2([0,1],C_0^2(D_j)))^4 \\ (u_2(t,x,y,z),v_2(t,x,y,z),w_2(t,x,y,z),p_2(t,x,y,z)) \in (C^2([1,2],C_0^2(D_j)))^4 \\ (u_3(t,x,y,z),v_3(t,x,y,z),\mathrm{w}_3(t,x,y,z),p_3(t,x,y,z)) \in (C^2([2,3],C_0^2(D_j)))^4 \\ (u_4(t,x,y,z),v_4(t,x,y,z),w_4(t,x,y,z),p_4(t,x,y,z)) \in (C^2([3,4],C_0^2(D_j)))^4 \\ \dots \end{cases}$$

belongs to $(C^1([0,\infty),C_0^2(D_j)))^4$ and it is a solution to the problem (3.2).

Note that

$$\mathrm{supp}(u^j), \quad \mathrm{supp}(v^j), \quad \mathrm{supp}(w^j), \quad \mathrm{supp}(p^j) \subset D_{jj} \subset D_j \quad \text{for} \quad \text{any} \quad j = 1,2,\dots.$$

and then

$$\left(D_{txyz}^{\alpha}u^j, D_{txyz}^{\alpha}v^j, D_{txyz}^{\alpha}w^j, D_{txyz}^{\alpha}p^j\right)|_{\partial D_j}$$

$$= \left(D_{txyz}^{\alpha}u^{j+1}, D_{txyz}^{\alpha}v^{j+1}, D_{txyz}^{\alpha}w^{j+1}, D_{txyz}^{\alpha}p^{j+1}\right)|_{\partial D_{j+1}}$$

$$= 0$$

for any $\alpha = (\alpha_0,\alpha_1,\alpha_2,\alpha_3)$, $\alpha_0,\alpha_1,\alpha_2,\alpha_3 \in \{0,1,\dots\}$. Also,

$$(u,v,w,p) = \begin{cases} (u^1,v^1,w^1,p^1) \in (C^2([0,\infty),C_0^2(D_1)))^4 \\ (u^2,v^2,w^2,p^2) \in (C^2([0,\infty),C_0^2(D_2)))^4 \\ (u^3,v^3,w^3,p^3) \in (C^2([0,\infty),C_0^2(D_3)))^4 \\ \dots \end{cases}$$

is a solution to the problem (1.1) which belongs to the space $(C^2([0,\infty)\times\mathbb{R}^3))^4$. Using the system (1.1) we have that $(u,v,w,p) \in (C^\infty([0,\infty)\times\mathbb{R}^3))^4$.

Therefore

$$\mathrm{supp}u, \mathrm{supp}v, \mathrm{supp}w, \mathrm{supp}p \subset D_{11} \cup D_{22} \cup \dots.$$

Also, for any $t \in [0,\infty)$, we have

$$\int_{\mathbb{R}^3} |u(t,x,y,z)|^2 dxdydz = \sum_{j=1}^{\infty} \int_{D_j} |u^j(t,x,y,z)|^2 dxdydz \le \sum_{j=1}^{\infty} \frac{1}{2^j} < \infty,$$

$$\int_{\mathbb{R}^3} |v(t,x,y,z)|^2 dxdydz = \sum_{j=1}^{\infty} \int_{D_j} |v^j(t,x,y,z)|^2 dxdydz \le \sum_{j=1}^{\infty} \frac{1}{2^j} < \infty,$$

$$\int_{\mathbb{R}^3} |w(t,x,y,z)|^2 dxdydz = \sum_{j=1}^{\infty} \int_{D_j} |w^j(t,x,y,z)|^2 dxdydz \le \sum_{j=1}^{\infty} \frac{1}{2^j} < \infty,$$

$$\int_{\mathbb{R}^3} |p(t,x,y,z)|^2 dxdydz = \sum_{j=1}^{\infty} \int_{D_j} |p^j(t,x,y,z)|^2 dxdydz \le \sum_{j=1}^{\infty} \frac{1}{2^j} < \infty.$$

This completes the proof.

**Remark 3.4** *We note that in mth step we have*

$$I_1^{mj}(u,v,w,p)=\int_{x_0}^{x}\int_{x_0}^{x_1}\int_{y_0}^{y}\int_{y_0}^{y_1}\int_{z_0}^{z}\int_{z_0}^{z_1}(u(t,\alpha,\beta,\gamma)$$

$$-u_{m-1}(m-1,\alpha,\beta,\gamma))d\gamma dz_1 d\beta dy_1 d\alpha dx_1$$

$$+\int_{m-1}^{t}\int_{x_0}^{x}\int_{y_0}^{y}\int_{y_0}^{y_1}\int_{z_0}^{z}\int_{z_0}^{z_1}u^2(s,\alpha,\beta,\gamma)d\gamma dz_1 d\beta dy_1 d\alpha ds$$

$$+\int_{m-1}^{t}\int_{x_0}^{x}\int_{x_0}^{x_1}\int_{y_0}^{y}\int_{z_0}^{z}\int_{z_0}^{z_1}u(s,\alpha,\beta,\gamma)v(s,\alpha,\beta,\gamma)d\gamma dz_1 d\beta d\alpha dx_1 ds$$

$$+\int_{m-1}^{t}\int_{x_0}^{x}\int_{x_0}^{x_1}\int_{y_0}^{y}\int_{y_0}^{y_1}\int_{z_0}^{z}u(s,\alpha,\beta,\gamma)w(s,\alpha,\beta,\gamma)d\gamma d\beta dy_1 d\alpha dx_1 ds$$

$$+\int_{m-1}^{t}\int_{x_0}^{x}\int_{y_0}^{y}\int_{y_0}^{y_1}\int_{z_0}^{z}\int_{z_0}^{z_1}p(s,\alpha,\beta,\gamma)d\gamma dz_1 d\beta dy_1 d\alpha ds$$

$$-\int_{m-1}^{t}\int_{y_0}^{y}\int_{y_0}^{y_1}\int_{z_0}^{z}\int_{z_0}^{z_1}u(s,x,\beta,\gamma)d\gamma dz_1 d\beta dy_1 ds$$

$$-\int_{m-1}^{t}\int_{x_0}^{x}\int_{x_0}^{x_1}\int_{z_0}^{z}\int_{z_0}^{z_1}u(s,\alpha,y,\gamma)d\gamma dz_1 d\alpha dx_1 ds$$

$$-\int_{m-1}^{t}\int_{x_0}^{x}\int_{x_0}^{x_1}\int_{y_0}^{y}\int_{y_0}^{y_1}u(s,\alpha,\beta,z)d\beta dy_1 d\alpha dx_1 ds,$$

$$I_2^{mj}(u,v,w,p)=\int_{x_0}^{x}\int_{x_0}^{x_1}\int_{y_0}^{y}\int_{y_0}^{y_1}\int_{z_0}^{z}\int_{z_0}^{z_1}(v(t,\alpha,\beta,\gamma)$$

$$-v_{m-1}(m-1,\alpha,\beta,\gamma))d\gamma dz_1 d\beta dy_1 d\alpha dx_1$$

$$+\int_{m-1}^{t}\int_{x_0}^{x}\int_{y_0}^{y}\int_{y_0}^{y_1}\int_{z_0}^{z}\int_{z_0}^{z_1}u(s,\alpha,\beta,\gamma)v(s,\alpha,\beta,\gamma)d\gamma dz_1 d\beta dy_1 d\alpha ds$$

$$+\int_{m-1}^{t}\int_{x_0}^{x}\int_{x_0}^{x_1}\int_{y_0}^{y}\int_{z_0}^{z}\int_{z_0}^{z_1}v^2(s,\alpha,\beta,\gamma)d\gamma dz_1 d\beta dy_1 d\alpha dx_1 ds$$

$$+\int_{m-1}^{t}\int_{x_0}^{\mathrm{x}}\int_{x_0}^{x_1}\int_{y_0}^{y}\int_{y_0}^{y_1}\int_{z_0}^{z}v(s,\alpha,\beta,\gamma)w(s,\alpha,\beta,\gamma)d\gamma d\beta dy_1 d\alpha dx_1 ds$$

$$+\int_{m-1}^{t}\int_{x_0}^{x}\int_{x_0}^{x_1}\int_{y_0}^{y}\int_{z_0}^{z}\int_{z_0}^{z_1}p(s,\alpha,\beta,\gamma)d\gamma dz_1 d\beta d\alpha dx_1 ds$$

$$-\int_{m-1}^{t}\int_{y_0}^{y}\int_{y_0}^{y_1}\int_{z_0}^{z}\int_{z_0}^{z_1}v(s,x,\beta,\gamma)d\gamma dz_1 d\beta dy_1 ds$$

$$-\int_{m-1}^{t}\int_{x_0}^{x}\int_{x_0}^{x_1}\int_{z_0}^{z}\int_{z_0}^{z_1}v(s,\alpha,y,\gamma)d\gamma dz_1 d\alpha dx_1 ds$$

$$-\int_{m-1}^{t}\int_{x_0}^{x}\int_{x_0}^{x_1}\int_{y_0}^{y}\int_{y_0}^{y_1} v(s,\alpha,\beta,z)\,d\beta\, dy_1\, d\alpha\, dx_1\, ds,$$

$$I_3^{mj}(u,v,w,p)=\int_{x_0}^{x}\int_{x_0}^{x_1}\int_{y_0}^{y}\int_{y_0}^{y_1}\int_{z_0}^{z}\int_{z_0}^{z_1}(w(t,\alpha,\beta,\gamma)$$

$$-w_{m-1}(m-1,\alpha,\beta,\gamma))\,d\gamma\, dz_1\, d\beta\, dy_1\, d\gamma\, dx_1$$

$$+\int_{m-1}^{t}\int_{x_0}^{x}\int_{y_0}^{y}\int_{y_0}^{y_1}\int_{z_0}^{z}\int_{z_0}^{z_1} u(s,\alpha,\beta,\gamma)w(s,\alpha,\beta,\gamma)\,d\gamma\, dz_1\, d\beta\, dy_1\, d\alpha\, ds$$

$$+\int_{m-1}^{t}\int_{x_0}^{x}\int_{x_0}^{x_1}\int_{y_0}^{y}\int_{z_0}^{z}\int_{z_0}^{z_1} v(s,\alpha,\beta,\gamma)w(s,\alpha,\beta,\gamma)\,d\gamma\, dz_1\, d\beta\, d\alpha\, dx_1\, ds$$

$$+\int_{m-1}^{t}\int_{x_0}^{x}\int_{x_0}^{x_1}\int_{y_0}^{y}\int_{y_0}^{y_1}\int_{z_0}^{z} w^2(s,\alpha,\beta,\gamma)\,d\gamma\, d\beta\, dy_1\, d\alpha\, dx_1\, ds$$

$$+\int_{m-1}^{t}\int_{x_0}^{x}\int_{x_0}^{x_1}\int_{y_0}^{y}\int_{y_0}^{y_1}\int_{z_0}^{z} p(s,\alpha,\beta,\gamma)\,d\gamma\, d\beta\, dy_1\, d\alpha\, dx_1\, ds$$

$$-\int_{m-1}^{t}\int_{y_0}^{y}\int_{y_0}^{y_1}\int_{z_0}^{z}\int_{z_0}^{z_1} w(s,x,\beta,\gamma)\,d\gamma\, dz_1\, d\beta\, dy_1\, ds$$

$$-\int_{m-1}^{t}\int_{x_0}^{x}\int_{x_0}^{x_1}\int_{z_0}^{z}\int_{z_0}^{z_1} w(s,\alpha,y,\gamma)\,d\gamma\, dz_1\, d\alpha\, dx_1\, ds$$

$$-\int_{m-1}^{t}\int_{x_0}^{x}\int_{x_0}^{x_1}\int_{y_0}^{y}\int_{y_0}^{y_1} w(s,\alpha,\beta,z)\,d\beta\, dy_1\, d\alpha\, dx_1\, ds,$$

$$I_4^{mj}(u,v,w,p)=\int_{0}^{t}\int_{x_0}^{x}\int_{y_0}^{y}\int_{y_0}^{y_1}\int_{z_0}^{z}\int_{z_0}^{z_1} u(s,\alpha,\beta,\gamma)\,d\gamma\, dz_1\, d\beta\, dy_1\, d\alpha\, ds$$

$$+\int_{m-1}^{t}\int_{x_0}^{x}\int_{x_0}^{x_1}\int_{y_0}^{y}\int_{z_0}^{z}\int_{z_0}^{z_1} v(s,\alpha,\beta,\gamma)\,d\gamma\, dz_1\, d\beta\, d\alpha\, dx_1\, ds$$

$$+\int_{m-1}^{t}\int_{x_0}^{x}\int_{x_0}^{x_1}\int_{y_0}^{y}\int_{y_0}^{y_1}\int_{z_0}^{z} w(s,\alpha,\beta,\gamma)\,d\gamma\, d\beta\, dy_1\, d\alpha\, dx_1\, ds.$$

## 4 Acknowledgements

We would like to thank Professor Simeon Reich, Department of Mathematics, Technion, Israel for his insightful comments.

# Existence and Smoothness of a Class of Burgers Equations

Svetlin G. Georgiev and Gal Davidi

Mon Jan 18 22:39:33 2021

In this paper we consider a class of Burgers equations. We propose a new method for investigation for existence of classical solutions.

## 1 Introduction

When solving a wide range of problems related to nonlinear acoustics, we may describe the nonlinear sound waves in fluids by using the Burgers model equation. This equation is named after Johannes Martinus Burgers, who published an equation of this form in his paper concerning turbulent phenomena modelling in 1948 (see [2]). However, this type of equation used in continuum mechanics was first presented in a paper published in ameteorology journal by H. Bateman in 1915 (see [3]). Thanks to the fact that this journal was not read by experts in continuum mechanics, the equation has become known as the Burgers equation. The possibilities of using this equation in nonlinear acoustics were probably discovered by J. Cole in 1949. Subsequently E. Hopf (in 1950) and J. Cole (in 1951) published the transformation independently in their papers (see [4, 5]). The transformation enables us to find the general analytic solution of the Burgers equation, as a result of which this equation plays a major role in nonlinear acoustics. The Burgers equation can be used for modelling both travelling and standing nonlinear plane waves [1, 6, 7, 8, 10, 15]. The equation is the simplest model equation that can describe the second order nonlinear effects connected with the propagation of high-amplitude (finite-amplitude waves) plane waves and, in addition, the dissipative effects in real fluids [9, 11, 12, 13, 14, 18]. There are several approximate solutions of the Burgers equation (see [16]). These solutions are always fixed to areas before and after the shock formation. For an area where the shock wave is forming no approximate solution has yet been found. It is therefore necessary to solve the Burgers equation numerically in this area (see [17, 19]). Numerical solutions themselves have difficulties with stability and accuracy.

Here we focus our attention on a class of Burgers equations and we will investigate it for existence of classical solutions. More precisely, consider the IVP

$$\begin{aligned} u_t + uu_x &= 0, \quad t > 0, \quad x \in \mathbb{R}, \\ u(0,x) &= \phi(x), \end{aligned} \tag{1.1}$$

where

$$\phi \in \bigcup_{j=-\infty}^{\infty} \mathcal{C}_0^1([j,j+1]), \quad \| \phi \|_{L^2(\mathbb{R})} \leq 1,$$

$$\max\left(\sup_{x\in\mathbb{R}}|\phi(x)|, \quad \sup_{x\in\mathbb{R}}|\phi'(x)|\right) < 1, \tag{1.2}$$

$$\max\left(\max_{x\in[j,j+1]}|\phi(x)|, \quad \max_{x\in[j,j+1]}|\phi'(x)|\right) \leq \frac{1}{2^{|j|+1}}, \quad j \in \mathbb{Z}.$$

Our main result is as follows.

**Theorem 1.1** *Assume (1.2). Then the IVP (1.1) has a solution* $u \in \mathcal{C}^1([0,\infty), \bigcup_{j=-\infty}^{\infty} \mathcal{C}_0^1([j,j+1]))$.

The paper is organized as follows. In the next section we will give some preliminaries which will be used for the proof of our main result. In Section 3 we give the proof of our main result. In Section 4 we make some conclusions.

## 2 Preliminaries

To prove our main result we will use a fixed point theorem for a sum two operators one of which is an expansive operator.

**Theorem 2.1 (, Corollary 2.4)** *Let $X$ be a nonempty closed convex subset of a Banach space $E$. Suppose that $T$ and $S$ map $X$ into $E$ such that*

1. $S$ is continuous and $S(X)$ resides in a compact subset of $E$.
2. $T: X \to E$ is expansive and onto.

Then there exists a point $x^* \in X$ such that

$$Sx^* + Tx^* = x^*.$$

**Definition 2.2** *Let $(X,d)$ be a metric space and $M$ be a subset of $X$. The mapping $T: M \to X$ is said to be expansive if there exists a constant $h > 1$ such that*

$$d(Tx,Ty) \geq h d(x,y)$$

for any $x, y \in M$.

## 3 Proof of the Main Result

Let $\varepsilon \in \left(0, \frac{1}{10}\right)$.

For $u \in \mathcal{C}^1([0,1], \mathcal{C}_0^1([0,1]))$, we define the operators

$$T_{11}u(t,x) = (1+\varepsilon)u(t,x),$$

$$S_{11}u(t,x) = -\varepsilon u(t,x)$$

$$+\varepsilon\left(\int_0^x u(t,z)dz - \int_0^x \phi(z)dz + \frac{1}{2}\int_0^t (u(\tau,x))^2 d\tau\right),$$

$t \in [0,1], \quad x \in [0,1]$. Consider the IVP

$$\begin{aligned} u_t + u u_x &= 0, \quad (t,x) \in (0,1] \times [0,1], \\ u(0,x) &= \phi(x), \quad x \in [0,1]. \end{aligned} \tag{3.1}$$

Note that every $u \in \mathcal{C}^1([0,1], \mathcal{C}_0^1([0,1]))$ which is a fixed point of $T_{11} + S_{11}$ is a solution of the

IVP (3.1). Really, we have

$$u(t,x)=(1+\varepsilon)u(t,x)-\varepsilon u(t,x)$$

$$+\varepsilon\left(\int_0^x u(t,z)dz-\int_0^x \phi(z)dz+\frac{1}{2}\int_0^t (u(\tau,x))^2 d\tau\right),$$

$(t,x)\in[0,1]\times[0,1]$, whereupon

$$\int_0^x u(t,z)dz-\int_0^x \phi(z)dz+\frac{1}{2}\int_0^t (u(\tau,x))^2 d\tau=0,\quad (t,x)\in[0,1]\times[0,1]. \tag{3.2}$$

We differentiate the last equation with respect to $x$ and we get

$$u(t,x)-\phi(x)+\int_0^t u(\tau,x)u_x(\tau,x)d\tau=0,\quad (t,x)\in[0,1]\times[0,1].$$

Now we differentiate the last equation with respect to $t$ and we obtain

$$u_t(t,x)+u(t,x)u_x(t,x)=0,\quad (t,x)\in[0,1]\times[0,1].$$

We put $t=0$ in (3.2). Then

$$\int_0^x u(0,z)dz-\int_0^x \phi(z)dz=0,\quad x\in[0,1],$$

which we differentiate with respect to $x$ and we get

$$u(0,x)=\phi(x),\quad x\in[0,1].$$

Let $\tilde{\tilde{\tilde{X}}}^1$ be the set of all equicontinuous families of the space $\mathcal{C}^1([0,1],\mathcal{C}_0^1([0,1]))$ with respect to the norm

$$\| f \|==\max\{\max_{(t,x)\in[0,1]\times[0,1]}|u(t,x)|,\quad \max_{(t,x)\in[0,1]\times[0,1]}|u_t(t,x)|$$

$$\max_{(t,x)\in[0,1]\times[0,1]}|u_x(t,x)|\},$$

$\tilde{\tilde{X}}^1=\tilde{\tilde{\tilde{X}}}^1\cup\{\phi\}$, $\tilde{X}^1=\overline{\tilde{\tilde{X}}^1}$, i.e., $\tilde{X}^1$ is the completion of $\tilde{\tilde{X}}^1$,

$$X^1=\{f\in\tilde{X}^1:\| f \|\leq\frac{1}{2}\},$$

$$Y^1=\{f\in\tilde{X}^1:\| f \|\leq\frac{1+\varepsilon}{2}\}.$$

Note that $X^1\subset Y^1$ and $Y^1$ is a compact subset of $\mathcal{C}^1([0,1],\mathcal{C}_0^1([0,1]))$. For any $x,y\in X^1$, we have

$$\| T_{11}x-T_{11}y \|=(1+\varepsilon)\| x-y \|.$$

Consequently $T_{11}:X^1\to\mathcal{C}^1([0,1],\mathcal{C}_0^1([0,1]))$ is a linear operator which is expansive with a constant $1+\varepsilon$. Let $v\in Y^1$. Then $u=\frac{v}{1+\varepsilon}\in X^1$ and $T_{11}u=v$. Therefore $T_{11}:X^1\to Y^1$ is onto. Now we will prove that $S_{11}:X^1\to X^1$. Let $u\in X^1$. We have

$$|\phi(x)|\leq\frac{1}{2},\quad x\in[0,1],$$

$$|u(t,x)|\leq\frac{1}{2},\quad |u_t(t,x)|\leq\frac{1}{2},\quad |u_x(t,x)|\leq\frac{1}{2},\quad (t,x)\in[0,1]\times[0,1],$$

and

$$|S_{11}u(t,x)|\leq|-\varepsilon u(t,x)$$

$$+\varepsilon\left(\int_0^x u(t,z)dz-\int_0^{\mathrm{x}} \phi(z)dz+\frac{1}{2}\int_0^t (u(t,x))^2 dx\right)|$$

$$\leq\varepsilon(|u(t,x)|+\int_0^x |u(t,z)|dz+\int_0^x |\phi(z)|dz$$

$$+\frac{1}{2}\int_0^t |u(\tau,x)|^2 d\tau)$$

$$\leq \varepsilon\left(\frac{1}{2}+\frac{1}{2}+\frac{1}{2}+\frac{1}{8}\right)$$

$$=\frac{13}{8}\varepsilon$$

$$\leq \frac{1}{2}, \quad (t,x)\in[0,1]\times[0,1],$$

and

$$\frac{\partial}{\partial t}(S_{11}u(t,x)) = \varepsilon\left(-u_t(t,x)+\int_0^x u_t(t,z)dz+\frac{1}{2}(u(t,x))^2\right),$$

$(t,x)\in[0,1]\times[0,1]$. So,

$$|\frac{\partial}{\partial t}(S_{11}u(t,x))| \leq \varepsilon\left(|u_t(t,x)|+\int_0^x |u_t(t,z)|dz+\frac{1}{2}(u(t,x))^2\right)$$

$$\leq \varepsilon\left(\frac{1}{2}+\frac{1}{2}+\frac{1}{8}\right)$$

$$=\frac{9}{8}\varepsilon$$

$$\leq \frac{1}{2}, \quad (t,x)\in[0,1]\times[0,1].$$

Next,

$$\frac{\partial}{\partial x}(S_{11}u(t,x)) = \varepsilon\left(-u_x(t,x)+u(t,x)-\phi(x)+\int_0^t u(\tau,x)u_x(\tau,x)d\tau\right),$$

$(t,x)\in[0,1]\times[0,1]$. So,

$$\left|\frac{\partial}{\partial x}(S_{11}u(t,x))\right| \leq \varepsilon\left(|u_x(t,x)|+|u(t,x)|+|\phi(x)|+\int_0^t |u(\tau,x)u_x(\tau,x)|d\tau\right)$$

$$\leq \varepsilon\left(\frac{1}{2}+\frac{1}{2}+\frac{1}{2}+\frac{1}{4}\right)$$

$$=\frac{7}{4}\varepsilon$$

$$\leq \frac{1}{2}, \quad (t,x)\in[0,1]\times[0,1].$$

Consequently $S_{11}: X^1 \to X^1$. Hence and Theorem 2.1, it follows that the IVP (3.1) has a solution $u^{11}\in\mathcal{C}^1([0,1],\mathcal{C}_0^1([0,1]))$. Now, for $u\in\mathcal{C}^1([0,1],\mathcal{C}_0^1([1,2]))$, we define the operators

$$T_{12}u(t,x) = (1+\varepsilon_1)u(t,x),$$

$$S_{12}u(t,x) = -\varepsilon_1 u(t,x)$$

$$+\varepsilon\left(\int_1^x u(t,z)dz-\int_1^x \phi(z)dz+\frac{1}{2}\int_0^t ((u(\tau,x))^2-(u^{11}(\tau,1))^2)d\tau\right),$$

$t\in[0,1], \quad x\in[1,2]$, where $\varepsilon_1\in\left(0,\frac{1}{10^2}\right)$. Consider the IVP

$$\begin{aligned} u_t+uu_x &= 0, \quad (t,x)\in[0,1]\times[1,2], \\ u(0,x) &= \phi(x), \quad x\in[1,2]. \end{aligned} \tag{3.3}$$

Note that every $u \in \mathcal{C}^1([0,1], \mathcal{C}_0^1([1,2]))$ which is a fixed point of $T_{12} + S_{12}$ is a solution of the IVP (3.3). Let $\tilde{\tilde{\tilde{X}}}^2$ be the set of all equicontinuous families of the space $\mathcal{C}^1([0,1], \mathcal{C}_0^1([1,2]))$ with respect to the norm

$$\| f \| == \max\{\max_{(t,x)\in[0,1]\times[1,2]} |u(t,x)|, \quad \max_{(t,x)\in[0,1]\times[1,2]} |u_t(t,x)|$$

$$\max_{(t,x)\in[0,1]\times[1,2]} |u_x(t,x)|\},$$

$\tilde{\tilde{X}}^2 = \tilde{\tilde{\tilde{X}}}^2 \cup \quad \{\phi, u^{11}\}$, $\tilde{X}^2 = \overline{\tilde{\tilde{X}}^2}$, i.e., $\tilde{X}^2$ is the completion of $\tilde{\tilde{X}}^2$,

$$X^2 = \{f \in \tilde{X}^2 : \| f \| \leq \frac{1}{2^2}\},$$

$$Y^2 = \{f \in \tilde{X}^2 : \| f \| \leq \frac{1+\varepsilon_1}{2^2}\}.$$

Note that $X^2 \subset Y^2$ and $Y^2$ is a compact subset of $\mathcal{C}^1([0,1], \mathcal{C}_0^1([1,2]))$. As in above we get a solution $u^{12} \in \mathcal{C}^1([0,1], \mathcal{C}_0^1([1,2]))$ of the IVP (3.3). Note that

$$\int_1^x u^{12}(t,z)dz - \int_1^x \phi(z)dz + \frac{1}{2}\int_0^t ((u^{12}(\tau,x))^2 - (u^{11}(\tau,1))^2)d\tau = 0,$$

$(t,x) \in [0,1] \times [0,1]$. Since $X^1 \subset \mathcal{C}^1([0,1], \mathcal{C}_0^1([0,1]))$ and $X^2 \subset \mathcal{C}^1([0,1], \mathcal{C}_0^1([1,2]))$, we have

$$u^{11}(t,1) = u_x^{11}(t,1) = 0, \quad t \in [0,1],$$

and

$$u^{12}(t,1) = u_x^{12}(t,1) = 0, \quad t \in [0,1].$$

Therefore

$$u^{11}(t,1) = u^{12}(t,1) = 0, \quad t \in [0,1], \tag{3.4}$$

and

$$u_x^{11}(t,1) = u_x^{12}(t,1) = 0, \quad t \in [0,1].$$

By (3.4), we find

$$u_t^{11}(t,1) = u_t^{12}(t,1) = 0, \quad t \in [0,1].$$

Thus, the function

$$v(t,x) = \begin{cases} u^{11}(t,x) & (t,x) \in [0,1] \times [0,1] \\ u^{12}(t,x) & (t,x) \in [0,1] \times [1,2] \end{cases}$$

is a $\mathcal{C}^1([0,1], \mathcal{C}_0^1([0,2]))$-solution of the IVP

$$u_t + uu_x = 0 \quad (t,x) \in [0,1] \times [0,2]$$

$$u(0,x) = \phi(x), \quad x \in [0,2].$$

Next, for $u \in \mathcal{C}^1([0,1], \mathcal{C}_0^1([2,3]))$, we define the operators

$$T_{13}u(t,x) = (1+\varepsilon_2)u(t,x),$$

$$S_{13}u(t,x) = -\varepsilon_2 u(t,x)$$

$$+\varepsilon_2(\int_2^x u(t,z)dz - \int_2^x \phi(z)dz$$

$$+\frac{1}{2}\int_0^t ((u(\tau,x))^2 - (u^{12}(\tau,2))^2)d\tau),$$

$t \in [0,1], \quad x \in [2,3]$, where $\varepsilon \in \left(0, \frac{1}{10^4}\right)$. Consider the IVP

$$u_t + uu_x = 0, \quad (t,x) \in [0,1] \times [2,3],$$
$$u(0,x) = \phi(x), \quad x \in [2,3]. \tag{3.5}$$

Note that every $u \in \mathcal{C}^1([0,1], \mathcal{C}_0^1([2,3]))$ which is a fixed point of $T_{13} + S_{13}$ is a solution of the IVP (3.5). Let $\tilde{\tilde{\tilde{X}}}^3$ be the set of all equicontinuous families of the space $\mathcal{C}^1([0,1], \mathcal{C}_0^1([2,3]))$ with respect to the norm

$$\| f \| == \max\{\max_{(t,x)\in[0,1]\times[2,3]} |u(t,x)|, \quad \max_{(t,x)\in[0,1]\times[2,3]} |u_t(t,x)|$$

$$\max_{(t,x)\in[0,1]\times[2,3]} |u_x(t,x)|\},$$

$\tilde{\tilde{X}}^3 = \tilde{\tilde{\tilde{X}}}^3 \cup \{\phi, u^{12}\}$, $\tilde{X}^3 = \overline{\tilde{\tilde{X}}^3}$, i.e., $\tilde{X}^3$ is the completion of $\tilde{\tilde{X}}^3$,

$$X^3 = \{f \in \tilde{X}^2 : \| f \| \leq \frac{1}{2^3}\},$$

$$Y^2 = \{f \in \tilde{X}^2 : \| f \| \leq \frac{1+\varepsilon_2}{2^3}\}.$$

Note that $X^3$ is a compact subset of $Y^3$ and $Y^3$ is a compact subset of $\mathcal{C}^1([0,1], \mathcal{C}_0^1([2,3]))$. As in above we get a solution $u^{13} \in \mathcal{C}^1([0,1], \mathcal{C}_0^1([2,3]))$ of the IVP (3.5). Note that

$$\int_2^x u^{13}(t,z)dz - \int_2^x \phi(z)dz + \frac{1}{2}\int_0^t ((u^{13}(\tau,x))^2 - (u^{12}(\tau,2))^2)d\tau = 0,$$

$(t,x) \in [0,1] \times [0,1]$. Since $X^2 \subset \mathcal{C}^1([0,1], \mathcal{C}_0^1([1,2]))$ and $X^3 \subset \mathcal{C}^1([0,1], \mathcal{C}_0^1([2,3]))$, we have

$$u^{12}(t,2) = u_x^{12}(t,2) = 0, \quad t \in [0,1],$$

and

$$u^{13}(t,2) = u_x^{13}(t,2) = 0, \quad t \in [0,1].$$

Therefore

$$u^{12}(t,2) = u^{13}(t,2) = 0, \quad t \in [0,1], \tag{3.6}$$

and

$$u_x^{12}(t,2) = u_x^{13}(t,2) = 0, \quad t \in [0,1].$$

By (3.6), we find

$$u_t^{12}(t,2) = u_t^{13}(t,2) = 0, \quad t \in [0,1].$$

Thus, the function

$$v(t,x) = \begin{cases} u^{11}(t,x) & (t,x) \in [0,1] \times [0,1] \\ u^{12}(t,x) & (t,x) \in [0,1] \times [1,2] \\ u^{13}(t,x) & (t,x) \in [0,1] \times [2,3] \end{cases}$$

is a $\mathcal{C}^1([0,1], \mathcal{C}_0^1([0,3]))$-solution of the IVP

$$u_t + uu_x = 0 \quad (t,x) \in [0,1] \times [0,3]$$

$$u(0,x) = \phi(x), \quad x \in [0,3].$$

Continuing in this manner, we get a $\mathcal{C}^1([0,1], \cup_{j=-\infty}^{\infty} \mathcal{C}_0^1([j,j+1]))$-solution $u^1$ of the IVP

$$u_t + uu_x = 0, \quad (t,x) \in [0,1] \times \mathbb{R},$$

$$u(0,x) = \phi(x), \quad x \in \mathbb{R},$$

for which we have

$$\max\left(\max_{(t,x)\in[0,1]\times[j,j+1]}|u^1(t,x)|, \quad \max_{(t,x)\in[0,1]\times[j,j+1]}|u^1_t(t,x)|,\right.$$

$$\left.\max_{(t,x)\in[0,1]\times[j,j+1]}|u^1_x(t,s)|\right) \leq \frac{1}{2^{|j|+1}}, \quad j\in\mathbb{Z}.$$

Next, for $u\in\mathcal{C}^1([1,2],\mathcal{C}^1_0([0,1]))$, we define the operators

$$T_{21}u(t,x) = (1+\varepsilon)u(t,x),$$

$$S_{21}u(t,x) = -\varepsilon u(t,x)$$

$$+\varepsilon\left(\int_0^x u(t,z)dz - \int_0^x u^1(1,z)dz + \frac{1}{2}\int_0^t (u(\tau,x))^2 d\tau\right),$$

$t\in[0,1]$, $\quad x\in[0,1]$, where $\varepsilon$ is as in Step 1. Consider the IVP

$$\begin{aligned} u_t + uu_x &= 0, \quad (t,x)\in[1,2]\times[0,1], \\ u(1,x) &= \phi(x), \quad x\in[0,1]. \end{aligned} \tag{3.7}$$

Note that every $u\in\mathcal{C}^1([1,2],\mathcal{C}^1_0([0,1]))$ which is a fixed point of $T_{21}+S_{21}$ is a solution of the IVP (3.7). Let $\tilde{\tilde{\tilde{X}}}^{21}$ be the set of all equicontinuous families of the space $\mathcal{C}^1([1,2],\mathcal{C}^1_0([0,1]))$ with respect to the norm

$$\| f \| == \max\{\max_{(t,x)\in[1,2]\times[0,1]}|u(t,x)|, \quad \max_{(t,x)\in[1,2]\times[0,1]}|u_t(t,x)|$$

$$\max_{(t,x)\in[1,2]\times[0,1]}|u_x(t,x)|\},$$

$\tilde{\tilde{X}}^{21} = \tilde{\tilde{\tilde{X}}}^{21}\cup \quad \{u^1\}$, $\tilde{X}^{21} = \overline{\tilde{\tilde{X}}^{21}}$, i.e., $\tilde{X}^{21}$ is the completion of $\tilde{\tilde{X}}^{21}$,

$$X^{21} = \{f\in\tilde{X}^{21}: \| f \| \leq \tfrac{1}{2}\},$$

$$Y^{21} = \{f\in\tilde{X}^{21}: \| f \| \leq \tfrac{1+\varepsilon}{2}\}.$$

Note that $X^{21}\subset Y^{21}$ and $Y^{21}$ is a compact subset of $\mathcal{C}^1([1,2],\mathcal{C}^1_0([0,1]))$. As in above we get a solution $u^{21}\in\mathcal{C}^1([1,2],\mathcal{C}^1_0([0,1]))$ of the IVP (3.7). Note that

$$\int_0^x u^{21}(t,z)dz - \int_0^x u^1(1,z)dz + \frac{1}{2}\int_1^t (u^{21}(\tau,x))^2 d\tau = 0, \quad (t,x)\in[1,2]\times[0,1].$$

We put $t=1$ in the last equation and we find

$$\int_0^x u^{21}(1,z)dz - \int_0^x u^1(1,z)dz = 0, \quad x\in[0,1].$$

We differentiate with respect to $x$, $x\in[0,1]$, the last equation and we obtain

$$u^{21}(1,x) = u^1(1,x), \quad x\in[0,1],$$

whereupon

$$u^{21}_x(1,x) = u^1_x(1,x), \quad x\in[0,1].$$

Since

$$0 = u^{21}_t(1,x) - u^{21}(1,x)u^{21}_x(1,x)$$

$$= u^{21}_t(1,x) - u^1(1,x)u^1_x(1,x), \quad x\in[0,1],$$

and

$$0 = u^1_t(1,x) - u^1(1,x)u^1_x(1,x), \quad x\in[0,1],$$

we conclude that

$$u_t^{21}(1,x) = u_t^1(1,x), \quad x \in [0,1].$$

And so on, we get a $\mathcal{C}^1([1,2], \cup_{j=-\infty}^{\infty} \mathcal{C}_0^1([j,j+1]))$-solution $u^2$ of the IVP

$$u_t + uu_x = 0, \quad (t,x) \in [1,2] \times \mathbb{R},$$

$$u(0,x) = u^1(1,x), \quad x \in \mathbb{R},$$

for which we have

$$|u^2|, \quad |u_t^2|, \quad |u_x^2| \leq \frac{1}{2^{|j|+1}} \quad on \quad [1,2] \times [j,j+1], \quad j \in \mathbb{Z}.$$

And so on. This completes the proof.

## 4 Conclusions

In this paper we investigate a class of Burgers equations for existence of classical solutions. We use a fixed point result for sum of two operators one of which is an expansive operator. We give a new integral representation of the solutions of the considered class of Burgers equations. The proposed approach in the paper is a new approach for investigations for existence of classical solutions. It can be used for investigations for existence of classical solutions of other classes partial differential equations.

# Classical Solutions for a Class of Compressible Navier-Stokes Equations

Svetlin G. Georgiev and Gal Davidi

Mon Jan 18 22:50:41 2021

In this article we investigate a class of compressible Navier-Stokes equations. We propose a new method for investigation for existence of classical solutions. This method gives new results.

## 1 Introduction

In this paper we consider the following IVP for a class of Navier-Stokes equations

$$\frac{\partial \rho}{\partial t} + \sum_{i=1}^{3} \frac{\partial}{\partial x_i}(\rho u_i) = 0$$

$$\frac{\partial}{\partial t}(\rho u_j) + \sum_{i=1}^{3} \frac{\partial}{\partial x_i}(\rho u_i u_j - P_{ij}) = 0, \tag{1.1}$$

$$\frac{\partial E}{\partial t} + \sum_{i=1}^{3} \left(u_i E - \sum_{j=1}^{3} u_j P_{ij} + q_i\right) = 0, \quad j \in \{1,2,3\}, \quad t > 0, \quad x \in \mathbb{R}^3,$$

$$\rho(0,x) = \rho_0(x), \quad E(0,x) = E_0(x), \quad u(0,x) = u_0(x), \quad x \in \mathbb{R}^3,$$

where $\rho$ is the mass density, $u$ is the velocity, $P_{ij}$, $i,j \in \{1,2,3\}$, is the stress tensor, $E$ is the energy density, $q$ is the heat flux and it is approximated by the Fourier law $q = -\kappa \nabla T$,

$$P_{ij} = -nk_0 T \delta_{ij} - \frac{2}{3}\mu \sum_{k=1}^{3} \frac{\partial u_k}{\partial x_k}\delta_{ij} + \mu\left(\frac{\partial u_i}{\partial x_j} + \frac{\partial u_j}{\partial x_i}\right), \quad i,j \in \{1,2,3\},$$

$$E = \frac{3}{2}nk_B T + \frac{1}{2}\rho u^2,$$

the density $n = \frac{\rho}{m}$, $m$ is the mass, and $T$ is the temperature, $\mu$ has units of $\rho\nu$, $\nu$ has units of $\frac{L^2}{t}$, $\delta_{ij}$, $i,j \in \{1,2,3\}$, is the Kronecker symbol,

$\rho_0, E_0 \in \mathcal{C}^2(\mathbb{R}^3)$, $u_0 \in (\mathcal{C}^2(\mathbb{R}^3))^3$,

$$|\rho_0| \le A, \quad |E_0| \le A, \quad |u_0| \le A,$$

where $A$ is a positive constant such that

$$\frac{3}{2}\frac{A}{m}k_B|T| + \frac{1}{2}A^3 \le A.$$

Here we investigate the system (1.1) for existence of classical solutions. Our main result is as follows.

**Theorem 1.1** *Suppose that $(H1)$ holds. Then the system (1.1) has at least one solution $(\rho, u, E)$ such that $\rho, E \in \mathcal{C}^1([0,\infty), \mathcal{C}^3(\mathbb{R}^3))$, $u \in (\mathcal{C}^1([0,\infty), \mathcal{C}^3(\mathbb{R}^3)))^3$.*

The Navier-Stokes model describes the variation of density at the interfaces between two phases, generally a liquid–vapor mixture(see [9-14] for more details regarding the physical illustrations). The existence of weak solutions for the compressible Navier-Stokes equations is an open problem in the case when viscosity coefficients depend on density. The existence of solutions of the compressible Navier-Stokes equations is obtained by Galerkin method(see [15-18] and references therein). The construction of approximation solutions to the Navier-Stokes equations is complicated. Antonelli and Spirito [19] have proved the global existence of weak solutions when $k < \nu$. In this paper we propose a new approach for investigation of existence of solutions for 3D NS equations.

The paper is organized as follows. In the next section we will give some preliminaries which will be used for the proof of our main result. In Section 3 we give the proof of our main result.

# 2 Preliminaries

To prove our main result we will use a fixed point theorem for a sum two operators one of which is an expansive operator.

**Definition 2.1** *Let $(X, d)$ be a metric space and $M$ be a subset of $X$. The mapping $T: M \to X$ is said to be expansive if there exists a constant $h > 1$ such that*

$$d(Tx, Ty) \geq hd(x, y)$$

for any $x, y \in M$.

**Theorem 2.2 (, Corollary 2.4)** *Let $X$ be a nonempty closed convex subset of a Banach space $E$. Suppose that $T$ and $S$ map $X$ into $E$ such that*

1. $S$ is continuous and $S(X)$ resides in a compact subset of $E$.
2. $T: X \to E$ is expansive and onto.

Then there exists a point $x^* \in X$ such that

$$Sx^* + Tx^* = x^*.$$

We introduce the following notations.

$$\int_0^{\hat{x}_l} = \int_0^{x_1} \cdots \int_0^{x_{l-1}} \int_0^{x_{l+1}} \cdots \int_0^{x_n},$$

$$dy_l = dy_1 \ldots dy_{l-1} dy_{l+1} \ldots dy_n,$$

$$z_{x_l} = (z_1, \ldots, z_{l-1}, x_l, z_{l+1}, \ldots, z_n), \quad l \in \{1, \ldots, n\},$$

$$dy = dy_1 \ldots dy_n.$$

# 3 Proof of the Main Result

Let $v \in \mathcal{C}([0,\infty) \times \mathbb{R}^3)$ be a function such that

$$\left| \int_0^t \int_0^x \int_0^y |v(s,z)| dz dy ds \right| \leq B,$$

$$\left|\int_0^x \int_0^y |v(t,z)| dz dy\right| \le B,$$

$$\left|\int_0^t \int_0^x \int_0^{\hat{y}_l} |v(s, z_{y_l})| dz_l dy ds\right| \le B,$$

$$\left|\int_0^t \int_0^{\hat{x}_l} \int_0^{\hat{y}_l} |v(s, z_{x_l})| dz_l dy_l ds\right| \le B,$$

$$\left|\int_0^t \int_0^x \int_0^y s|v(s,z)| dz dy ds\right| \le B,$$

$$\left|\int_0^x \int_0^y t|v(t,z)| dz dy\right| \le B,$$

$$\left|\int_0^t \int_0^x \int_0^{\hat{y}_l} s|v(s, z_{y_l})| dz_l dy ds\right| \le B,$$

$$\left|\int_0^t \int_0^{\hat{x}_l} \int_0^{\hat{y}_l} s|v(s, z_{x_l})| dz_l dy_l ds\right| \le B, \quad l \in \{1,2,3\},$$

$(t,x) \in [0,\infty) \times \mathbb{R}^3$, where $B$ is a positive constant. We choose $\varepsilon > 0$ such that

$$\max\{\varepsilon A + \varepsilon(2AB + 6BA^2),$$

$$\varepsilon A + \varepsilon\left(2BA^2 + 3B\left(A^3 + \frac{A}{m}|k_B||T| + \frac{2}{3}|\mu| + A + 2|\mu|A\right)\right),$$

$$\varepsilon A + \varepsilon(2AB + 3B\left(A^2 + 3a\left(\frac{A}{m}|k_B||T| + \frac{2}{3}|\mu| + A + 2|\mu|A\right) + |q|\right))\}$$

$$\le A.$$

Let $\tilde{\tilde{\tilde{X}}}$ be the set of all equicontinuous families of functions in the space $(\mathcal{C}^1([0,\infty), \mathcal{C}^2(\mathbb{R}^3)))^5$ with respect to the norm

$$\| w \| = \max_{j\in\{1,\dots,5\}}\{ \sup_{(t,x)\in[0,\infty)\times\mathbb{R}^3} |w_j(t,x)|,$$

$$\sup_{(t,x)\in[0,\infty)\times\mathbb{R}^3} \left|\frac{\partial}{\partial t} w_j(t,x)\right|,$$

$$\sup_{(t,x)\in[0,\infty)\times\mathbb{R}^3} \left|\frac{\partial}{\partial x_l} w_j(t,x)\right|,$$

$$\sup_{(t,x)\in[0,\infty)\times\mathbb{R}^3} \left|\frac{\partial^2}{\partial x_l^2} w_j(t,x)\right|, \quad l \in \{1,2,3\}\},$$

$\tilde{\tilde{X}} = \tilde{\tilde{\tilde{X}}} \cup \quad \{\rho_0, u_{01}, u_{02}, u_{03}, E_0\}$, $\tilde{X} = \overline{\tilde{\tilde{X}}}$, i.e., $\tilde{X}$ is the cvompletion of $\tilde{\tilde{X}}$,

$$X = \{w \in \tilde{X} : \| w \| \le B\},$$

$$Y = \{w \in \tilde{X} : \| w \| \le (1+\varepsilon)B\}.$$

Note that $X$ and $Y$ are compact sets in $(\mathcal{C}^1([0,\infty), \mathcal{C}^2(\mathbb{R}^3)))^5$. For $(\rho, u, E) \in (\mathcal{C}^1([0,\infty), \mathcal{C}^2(\mathbb{R}^3)))^5$, $u = (u_1, u_2, u_3)$, define the operators.

$$S_0(\rho,u,E)(t,x) = -\varepsilon\rho(t,x) + \varepsilon(\int_0^t \int_0^x \int_0^y v(s,z)(\rho(s,z)-\rho_0(z))dzdyds$$

$$-\sum_{i=1}^3 \int_0^t \int_0^x \int_0^y v(t_1,z)\int_0^{t_1} \frac{\partial}{\partial x_i}(\rho u_i)(s,z)dsdzdydt_1),$$

$$S_j(\rho,u,E)(t,x) = -\varepsilon u_j(t,x) + \varepsilon(\int_0^t \int_0^x \int_0^y v(s,z)(\rho(s,z)u_j(s,z)-\rho_0(z)u_{0j}(z))dzdyds$$

$$+\sum_{i=1}^3 \int_0^t \int_0^x \int_0^y v(t_1,z)\int_0^{t_1} \frac{\partial}{\partial x_i}(\rho u_i u_j - P_{ij})(s,z)dsdzdydt_1), \quad j \in \{1,2,3\},$$

$$S_4(\rho,u,E)(t,x) = -\varepsilon E(t,x) + \varepsilon(\int_0^t \int_0^x \int_0^y v(t_1,z)(E(t_1,z)-E_0(z))dzdydt_1$$

$$+\sum_{i=1}^3 \int_0^t \int_0^x \int_0^y v(t_1,z)\int_0^{t_1} (u_i(s,z)E(s,z) - \sum_{j=1}^3 u_j(s,z)P_{ij}(s,z) + q_i)dsdzdydt_1),$$

$S = (S_0, S_1, S_2, S_3, S_4)$,
$(t,x) \in [0,\infty) \times \mathbb{R}^3$, and

$$T_0(\rho,u,E)(t,x) = (1+\varepsilon)\rho(t,x),$$

$$T_j(\rho,u,E)(t,x) = (1+\varepsilon)u_j(t,x), \quad j \in \{1,2,3\},$$

$$T_4(\rho,u,E)(t,x) = (1+\varepsilon)E(t,x), \quad (t,x) \in [0,\infty) \times \mathbb{R}^3,$$

$$T = (T_0, T_1, T_2, T_3, T_4).$$

Note that any fixed point $(\rho,u,E) \in (\mathcal{C}^1([0,\infty),\mathcal{C}^2(\mathbb{R}^3)))^5$ of the operator $T+S$ is a solutions of the system (1.1). Really, let $(\rho,u,E) \in (\mathcal{C}^1([0,\infty),\mathcal{C}^2(\mathbb{R}^3)))^5$ is a fixed point of the operator $T+S$. Then

$$0 = \int_0^t \int_0^x \int_0^y v(s,z)(\rho(s,z)-\rho_0(z))dzdyds$$

$$-\sum_{i=1}^3 \int_0^t \int_0^x \int_0^y v(t_1,z)\int_0^{t_1} \frac{\partial}{\partial x_i}(\rho u_i)(s,z)dsdzdydt_1,$$

$$0 = \int_0^t \int_0^x \int_0^y v(s,z)(\rho(s,z)u_j(s,z)-\rho_0(z)u_{0j}(z))dzdyds$$

$$+\sum_{i=1}^3 \int_0^t \int_0^x \int_0^y v(t_1,z)\int_0^{t_1} \frac{\partial}{\partial x_i}(\rho u_i u_j - P_{ij})(s,z)dsdzdydt_1, \quad j \in \{1,2,3\},$$

$$0 = \int_0^t \int_0^x \int_0^y v(t_1,z)(E(t_1,z)-E_0(z))dzdydt_1$$

$$+\sum_{i=1}^3 \int_0^t \int_0^x \int_0^y v(t_1,z)\int_0^{t_1} (u_i(s,z)E(s,z) - \sum_{j=1}^3 u_j(s,z)P_{ij}(s,z) + q_i)dsdzdydt_1,$$

$(t,x) \in [0,\infty) \times \mathbb{R}^3$. We differentiate the last system once in $t$, twice in $x_1$, ..., $x_3$ and we get

$$0 = v(tx)(\rho(t,x) - \rho_0(x))$$

$$-\sum_{i=1}^{3} v(t,x) \int_0^t \frac{\partial}{\partial x_i}(\rho u_i)(s,x),$$

$$0 = v(t,x)(\rho(t,x)u_j(t,x) - \rho_0(x)u_{0j}(x))$$

$$+\sum_{i=1}^{3} v(t,x) \int_0^t \frac{\partial}{\partial x_i}(\rho u_i u_j - P_{ij})(s,x)ds, \quad j \in \{1,2,3\},$$

$$0 = v(t,x)(E(t,x) - E_0(x))$$

$$+\sum_{i=1}^{3} v(t,x) \int_0^t (u_i(s,x)E(s,x) - \sum_{j=1}^{3} u_j(s,x)P_{ij}(s,x) + q_i)ds,$$

$(t,x) \in [0,\infty) \times \mathbb{R}^3$, whereupon

$$0 = \rho(t,x) - \rho_0(x)$$

$$-\sum_{i=1}^{3} \int_0^t \frac{\partial}{\partial x_i}(\rho u_i)(s,x),$$

$$0 = \rho(t,x)u_j(t,x) - \rho_0(x)u_{0j}(x)$$

$$+\sum_{i=1}^{3} \int_0^t \frac{\partial}{\partial x_i}(\rho u_i u_j - P_{ij})(s,x)ds, \quad j \in \{1,2,3\}, \tag{3.1}$$

$$0 = E(t,x) - E_0(x)$$

$$+\sum_{i=1}^{3} \int_0^t (u_i(s,x)E(s,x) - \sum_{j=1}^{3} u_j(s,x)P_{ij}(s,x) + q_i)ds,$$

$(t,x) \in [0,\infty) \times \mathbb{R}^3$. Now we differentiate the system (3.1) once in $t$, twice in $x_1$, $x_2$, $x_3$ and we get the first five equations of the system (1.1). We put $t = 0$ in (3.1) and we get

$$0 = \rho(0,x) - \rho_0(x),$$

$$0 = \rho(0,x)u_j(0,x) - \rho_0(x)u_{0j}(x), \quad j \in \{1,2,3\},$$

$$0 = E(0,x) - E_0(x), \quad x \in \mathbb{R}^3,$$

from where

$$\rho(0,x) = \rho_0(x), \quad u(0,x) = u_0(x), \quad E(0,x) = E_0(x), \quad x \in \mathbb{R}^3.$$

Observe that $T: X \to Y$ is continuous and if $(\rho, u, E) \in Y$, then $\frac{1}{1+\varepsilon}(\rho, u, E) \in X$ and $T\left(\frac{1}{1+\varepsilon}(\rho, u, E)\right) = (\rho, u, E)$, i.e., $T: X \to Y$ is onto. By the choice of $\varepsilon$, we get for $(\rho, u, E) \in X$

$$\| S(\rho,u,E) \| \leq \max\{\varepsilon A+\varepsilon(2AB+6BA^2),$$

$$\varepsilon A+\varepsilon\left(2BA^2+3B\left(A^3+\frac{A}{m}|k_B||T|+\frac{2}{3}|\mu|+A+2|\mu|A\right)\right),$$

$$\varepsilon A+\varepsilon(2AB+3B\left(A^2+3a\left(\frac{A}{m}|k_B||T|+\frac{2}{3}|\mu|+A+2|\mu|A\right)+|q|\right))\}$$

$$\leq A.$$

Hence and Theorem 2.2, we conclude that $T+S$ has a fixed point in $X$.

# Classical Solutions for a Class of Euler Equations of Ideal Gas Dynamics

Svetlin G. Georgiev and Gal Davidi

Mon Jan 18 22:41:14 2021

In this paper we propose a new approach for investigation for existence of classical solutions for a class of Euler equations of ideal gas dynamics.

## 1 Introduction

The Euler equations are a class of quasilinear hyperbolic equations governing adiabatic and inviscid flow. The equations represent Cauchy equations of conservation of mass (continuity), and balance of momentum and energy, and we can consider as a particular case of Navier–Stokes equations with zero viscosity and zero thermal conductivity. In fact, Euler equations can be obtained by linearization of some more precise continuity equations like Navier–Stokes equations in a local equilibrium state given by a Maxwellian. The Euler equations can be applied to incompressible and to compressible flow.

The short time existence of a smooth solution for the 3d incompressible Euler equations has been obtained by Lichtenstein [7]. Existence, uniqueness, and propagation of the regularity of the initial data is investigated in [1]. The Beale–Kato–Majda result has been improved by Kozono [6]. Constantin, Fefferman and Majda [4] show that the variations in the direction of the vorticity may produce singularities. Additional relevant results can be found in Majda and Bertozzi [8] and Constantin [3]. In [2], [5] and references therein is proved blow up for smooth solutions with infinite energy of the 3d Euler equations. In [9] is proved existence, uniqueness, and global regularity for all solutions in 2d case.

In this paper we investigate the following Euler equations of ideal gas dynamics.

$$
\begin{aligned}
&\frac{\partial\rho}{\partial t}+\sum_{i=1}^{n}\frac{\partial}{\partial x_i}(\rho u_i)=0\\
&\frac{\partial}{\partial t}(\rho u_j)+\sum_{i=1}^{n}\frac{\partial}{\partial x_i}\left(\rho u_i u_j+P\delta_{ij}\right)=0,\quad j\in\{1,\dots,n\},\\
&\frac{\partial E}{\partial t}+\sum_{i=1}^{n}\frac{\partial}{\partial x_i}(u_i(E+\mathrm{P}))=0,\quad t>0,\quad x\in\mathbb{R}_+^n,\\
&\rho(0,x)=\rho_0(x)\\
&u_i(0,x)=u_{i0}(x),\quad i\in\{1,\dots,n\},\\
&E(0,x)=E_0(x),\quad x\in\mathbb{R}_+^n,
\end{aligned}
\tag{1.1}
$$

where $\mathbb{R}^n_+ = \{x = (x_1, \dots, x_n): x_i \geq 0, \quad i \in \{1, \dots, n\}\}$, $\gamma > 1$, $\delta_{ij} = 0$ if $i \neq j$, $\delta_{ij} = 1$ if $i = j$, $i, j \in \{1, \dots, n\}$, $P = (\gamma - 1)\left(E - \frac{1}{2}\rho u^2\right)$, $\rho_0, u_{i0}, E_0 \in \mathcal{C}^1(\mathbb{R}^n_+)$, $i \in \{1, \dots, n\}$, $u^2 = \sum_{i=1}^n u_i^2$.

Our main result is as follows.

**Theorem 1.1** *Let* $\rho_0, u_{i0}, E_0 \in \mathcal{C}^1(\mathbb{R}^n_+)$*,* $i \in \{1, \dots, n\}$*, and*

$$|\rho_0(x)| \leq B, \quad |u_{i0}(x)| \leq B, \quad |E_0(x)| \leq B,$$

$$\left|\frac{\partial}{\partial x_j}\rho_0(x)\right| \leq B, \quad \left|\frac{\partial}{\partial x_j}u_{i0}(x)\right| \leq B, \quad \left|\frac{\partial}{\partial x_j}E_0(x)\right| \leq B, \quad i, j \in \{1, \dots, n\},$$

$x \in \mathbb{R}^n_+$ for some positive constant. Then the IVP (1.1) has a solution $\rho, u_i, E \in \mathcal{C}^1([0, \infty) \times \mathbb{R}^n_+)$, $i \in \{1, \dots, n\}$.

The paper is organized as follows. In the next section we give some preliminaries used for the proof of our main result. In Section 3 we prove our main result. A conclusion is given in Section 4.

## 2 Preliminaries

To prove our main result we will use a fixed point theorem for a sum two operators one of which is an expansive operator.

**Theorem 2.1 (, Corollary 2.4)** *Let $X$ be a nonempty closed convex subset of a Banach space $E$. Suppose that $T$ and $S$ map $X$ into $E$ such that*

1. $S$ is continuous and $S(X)$ resides in a compact subset of $E$.
2. $T: X \to E$ is expansive and onto.

Then there exists a point $x^* \in X$ such that

$$Sx^* + Tx^* = x^*.$$

**Definition 2.2 ()** *Let $(X, d)$ be a metric space and $M$ be a subset of $X$. The mapping $T: M \to X$ is said to be expansive if there exists a constant $h > 1$ such that*

$$d(Tx, \mathrm{T}y) \geq hd(x, y)$$

for any $x, y \in M$.

We introduce the following notations.

$$\int_0^{\overline{y_i}} = \int_0^{y_1} \dots \int_0^{y_{i-1}} \int_0^{y_{i+1}} \dots \int_0^{y_n},$$

$$dz_{y_i} = dz_1 \dots dz_{i-1} dz_{i+1} \dots dz_n,$$

$$\tilde{z}_{y_i} = (z_1, \dots, z_{i-1}, y_i, z_{i+1}, \dots, z_n),$$

$$dy = dy_1 \dots dy_n.$$

Let $\phi \in \mathcal{C}^1([0, \infty) \times \mathbb{R}^n_+)$ be a positive continuously-differentiable function such that

$$\phi(t, x) \leq A,$$

$$\int_0^t \int_0^x \phi(\tau, y) \prod_{j=1}^n y_j dy d\tau \leq A,$$

$$\int_0^t \int_0^x \tau\phi(\tau, y) \prod_{j=1}^n y_j dy d\tau \leq A,$$

$$\int_0^x \phi(t,y) \prod_{j=1}^n y_j dy d\tau \le A,$$

$$t \int_0^x \phi(t,y) \prod_{j=1}^n y_j dy d\tau \le A,$$

$$x_k \int_0^t \int_0^{\overline{x_k}} \phi(\tau, \tilde{y}_{x_k}) \prod_{j=1,\ \ j\neq i}^n y_j dy_{x_k} dt \le A,$$

$$x_k \int_0^t \int_0^{\overline{x_k}} \phi(\tau, \tilde{y}_{x_k}) \tau \prod_{j=1,\ \ j\neq i}^n y_j dy_{x_k} dt \le A, \quad i,k \in \{1, \dots, n\},$$

$(t,x) \in [0,\infty) \times \mathbb{R}^n_+$, for some positive constant $A$. Such function $\phi$ exists. Indeed, let $\phi(t,x) = \frac{1}{(1+t)^{10}(1+x_1)^{10}(1+x_2)^{10} \dots (1+x_n)^{10}}$, $(t,x) \in [0,\infty) \times \mathbb{R}^n_+$.

## 3 Prof of the Main Result

In this section we will prove our main result. For this aim, we give an integral representation of the solution of the IVP (1.1). Then we construct two operators and using Theorem 2.1, we show that the sum of both operators has a fixed point which is a solution of the IVP (1.1). In the end of the section we give an example.

Let $\tilde{\tilde{\tilde{X}}}$ be the set of all equicontinuous families of the space $(\mathcal{C}^1([0,\infty) \times \mathbb{R}^n_+))^{n+2}$ with respect to the norm

$$\| v \| = \max\{ \sup_{x\in[0,\infty)\times\mathbb{R}^n_+} |v(t,x)|, \quad \sup_{x\in[0,\infty)\times\mathbb{R}^n_+} \left|\frac{\partial}{\partial t} v(t,x)\right|,$$

$$\sup_{x\in[0,\infty)\times\mathbb{R}^n_+} \left|\frac{\partial}{\partial x_j} v(t,x)\right|, \quad j \in \{1, \dots, n\}\}.$$

Let also,

$$\tilde{\tilde{X}} = \tilde{\tilde{\tilde{X}}} \cup \{u_{i0}, E_0, \rho_0, \quad i \in \{1, \dots, n\}\},$$

$$\tilde{X} = \overline{\tilde{\tilde{X}}},$$

i.e., $\tilde{X}$ is the completion of $\tilde{\tilde{X}}$,

$$X = \{v \in \tilde{X} \colon \| v \| \le B\},$$

$$Y = \{v \in \tilde{X} \colon \| v \| \le (1+\varepsilon)B\}.$$

Note that $X$ is a compact subset of $Y$. Observe that, if $(\rho, u, E) \in X$, then

$$|P| \le (\gamma - 1)\left(|E| + \frac{1}{2}|\rho| \sum_{j=1}^n u_j^2\right)$$

$$\le (\gamma - 1)\left(B + \frac{1}{2} n B^3\right),$$

$$\frac{\partial}{\partial t} P = (\gamma - 1)\left(\frac{\partial E}{\partial t} - \frac{1}{2} u^2 \frac{\partial}{\partial t}\rho - \rho \sum_{j=1}^n u_j \frac{\partial u_j}{\partial t}\right),$$

$$\left|\frac{\partial}{\partial t} P\right| \le (\gamma - 1)\left(\left|\frac{\partial E}{\partial t}\right| + \left|\frac{1}{2} u^2 \frac{\partial}{\partial t}\rho\right| + \rho \sum_{j=1}^n |u_j| \left|\frac{\partial u_j}{\partial t}\right|\right)$$

$$\le (\gamma - 1)\left(B + \frac{1}{2} n B^3 + n B^3\right)$$

$$= (\gamma - 1)\left(B + \frac{3}{2}nB^3\right),$$

$$\left|\frac{\partial}{\partial x_i}P\right| \leq (\gamma - 1)\left(\left|\frac{\partial E}{\partial x_i}\right| + \left|\frac{1}{2}u^2\frac{\partial}{\partial x_i}\rho\right| + \rho\sum_{j=1}^{n}|u_j|\left|\frac{\partial u_j}{\partial x_i}\right|\right)$$

$$\leq (\gamma - 1)\left(B + \frac{1}{2}nB^3 + nB^3\right)$$

$$= (\gamma - 1)\left(B + \frac{3}{2}nB^3\right) \quad i \in \{1, \dots, n\},$$

on $[0, \infty) \times \mathbb{R}_+^n$. Consequently, for $(\rho, u, E) \in X$, we have

$$\| P \| \leq (\gamma - 1)\left(B + \frac{3}{2}nB^3\right).$$

Now we choose $\varepsilon > 0$ small enough so that

$$\varepsilon(B + 2BA + nB^2A) \leq B,$$

$$\varepsilon\left(B + 2B^2A + nB^3A + A(\gamma - 1)n\left(B + \frac{3}{2}nB^3\right)\right) \leq B,$$

$$\varepsilon\left(B + 2BA + AB^2n + ABn(\gamma - 1)\left(B + \frac{3}{2}nB^3\right)\right) \leq B.$$

For $(\rho, u, E) \in X$, define the operators

$$T_0(\rho, u, E)(t, x) = (1 + \varepsilon)\rho(t, x),$$

$$S_0(\rho, u, E)(t, x) = -\varepsilon\rho(t, x)$$

$$+\varepsilon(\int_0^t \int_0^x \phi(\tau, y)\int_0^y (\rho(\tau, z) - \rho_0(z))dzdyd\tau$$

$$+\sum_{i=1}^{n}\int_0^t \int_0^x \phi(\tau, y)\int_0^\tau \int_0^{\overline{y_i}} \rho(s, \tilde{z}_{y_i})u_i(s, \tilde{z}_{y_i})dz_{y_i}dsdyd\tau),$$

$$T_j(\rho, u, E)(t, x) = (1 + \varepsilon)u_j(t, x),$$

$$S_j(\rho, u, E)(t, x) = -\varepsilon u_j(t, x)$$

$$+\varepsilon(\int_0^t \int_0^x \phi(\tau, y)\int_0^y (\rho(\tau, z)u_j(\tau, z) - \rho_0(z)u_{j0}(z))dzdyd\tau$$

$$+\sum_{i=1}^{n}\int_0^t \int_0^x \phi(\tau, y)\int_0^\tau \int_0^{\overline{y_i}} (\rho(s, \tilde{z}_{y_i})u_i(s, \tilde{z}_{y_i}) + P(s, \tilde{z}_{y_i})\delta_{ij})dz_{y_i}dsdyd\tau)$$

$$j \in \{1, \dots, n\},$$

$$T_{n+1}(\rho,u,E)(t,x)=(1+\varepsilon)E(t,x),$$

$$S_{n+1}(\rho,u,E)(t,x)=-\varepsilon E(t,x)$$

$$+\varepsilon\Big(\int_0^t\int_0^x\phi(\tau,y)\int_0^y\left(E(\tau,z)-E_0(z)\right)dzdyd\tau$$

$$+\sum_{i=1}^n\int_0^t\int_0^x\phi(\tau,y)\int_0^\tau\int_0^{\overline{y_i}}u_i(s,\tilde{z}_{y_i})(E+P)(s,\tilde{z}_{y_i})dz_{y_i}dsdyd\tau\Big),$$

$$T(\rho,u,E)(t,x)=(T_0(\rho,u,E)(t,x),T_1(\rho,u,E)(t,x),\dots,T_{n+1}(\rho,u,E)(t,x)),$$

$$S(\rho,u,E)(t,x)=(S_0(\rho,u,E)(t,x),S_1(\rho,u,E)(t,x),\dots,S_{n+1}(\rho,u,E)(t,x)).$$

Note that if $(\rho,u,E)\in X$ is a fixed point of the operator $T+S$, then it satisfies (1.1). Really, let $(\rho,u,E)\in X$. Then

$$0=\int_0^t\int_0^x\phi(\tau,y)\int_0^y\left(\rho(\tau,z)-\rho_0(z)\right)dzdyd\tau$$

$$+\sum_{i=1}^n\int_0^t\int_0^x\phi(\tau,y)\int_0^\tau\int_0^{\overline{y_i}}\rho(s,\tilde{z}_{y_i})u_i(s,\tilde{z}_{y_i})dz_{y_i}dsdyd\tau,$$

$$0=\int_0^t\int_0^x\phi(\tau,y)\int_0^y\left(\rho(\tau,z)u_j(\tau,z)-\rho_0(z)u_{j0}(z)\right)dzdyd\tau$$

$$+\sum_{i=1}^n\int_0^t\int_0^x\phi(\tau,y)\int_0^\tau\int_0^{\overline{y_i}}\left(\rho(s,\tilde{z}_{y_i})u_i(s,\tilde{z}_{y_i})+P(s,\tilde{z}_{y_i})\delta_{ij}\right)dz_{y_i}dsdyd\tau$$

$$j\in\{1,\dots,n\},$$

$$0=\int_0^t\int_0^x\phi(\tau,y)\int_0^y\left(E(\tau,z)-E_0(z)\right)dzdyd\tau$$

$$+\sum_{i=1}^n\int_0^t\int_0^x\phi(\tau,y)\int_0^\tau\int_0^{\overline{y_i}}u_i(s,\tilde{z}_{y_i})(E+P)(s,\tilde{z}_{y_i})dz_{y_i}dsdyd\tau.$$

We differentiate the last system with respect to $t$ and $x_i$, $i\in\{1,\dots,n\}$, and we get

$$0=\phi(t,x)\int_0^x\left(\rho(t,z)-\rho_0(z)\right)dz$$

$$+\sum_{i=1}^n\phi(t,x)\int_0^t\int_0^{\overline{x_i}}\rho(s,\tilde{z}_{x_i})u_i(s,\tilde{z}_{x_i})dz_{x_i}ds,$$

$$0=\phi(t,x)\int_0^x\left(\rho(t,z)u_j(t,z)-\rho_0(z)u_{j0}(z)\right)dz$$

$$+\sum_{i=1}^n\phi(t,x)\int_0^t\int_0^{\overline{x_i}}\left(\rho(s,\tilde{z}_{x_i})u_i(s,\tilde{z}_{x_i})+P(s,\tilde{z}_{x_i})\delta_{\mathrm{i}j}\right)dz_{x_i}ds$$

$j\in\{1,\dots,n\}$,

$$0=\phi(t,x)\int_0^x\left(E(t,z)-E_0(z)\right)dz$$

$$+\sum_{i=1}^n\phi(t,x)\int_0^t\int_0^{\overline{x_i}}u_i(s,\tilde{z}_{x_i})(E+P)(s,\tilde{z}_{x_i})dz_{x_i}ds,$$

or

$$0=\int_0^x\left(\rho(t,z)-\rho_0(z)\right)dz$$

$$+\sum_{i=1}^n\int_0^t\int_0^{\overline{x_i}}\rho(s,\tilde{z}_{x_i})u_i(s,\tilde{z}_{x_i})dz_{x_i}ds,$$

$$0=\int_0^x\left(\rho(t,z)u_j(t,z)-\rho_0(z)u_{j0}(z)\right)dz$$

$$+\sum_{i=1}^n\int_0^t\int_0^{\overline{x_i}}\left(\rho(s,\tilde{z}_{x_i})u_i(s,\tilde{z}_{x_i})+P(s,\tilde{z}_{x_i})\delta_{ij}\right)dz_{x_i}ds \qquad (3.1)$$

$j\in\{1,\dots,n\}$,

$$0=\int_0^x\left(E(t,z)-E_0(z)\right)dz$$

$$+\sum_{i=1}^n\int_0^t\int_0^{\overline{x_i}}u_i(s,\tilde{z}_{x_i})(E+P)(s,\tilde{z}_{x_i})dz_{x_i}ds.$$

Now we differentiate the last system with respect to $t$, $x_1$, …, $x_n$ and we get the first three equations of the system (1.1). Now we put $t=0$ in the system (3.1), and we get

$$0=\int_0^x\left(\rho(0,z)-\rho_0(z)\right)dz$$

$$0=\int_0^x\left(\rho(0,z)u_j(0,z)-\rho_0(z)u_{j0}(z)\right)dz$$

$$0=\int_0^x\left(E(0,z)-E_0(z)\right)dz.$$

We differentiate the last system with respect to $x_1$, $x_2$, …, $x_n$ and we obtain the last three equations of the system (1.1). Observe that, for $(\rho,u,E)\in X$, we have

$$\parallel S_0(\rho,u,E)\parallel\leq\varepsilon(B+2BA+nB^2A)$$

$$\leq B,$$

$$\| S_j(\rho,u,E) \| \leq \varepsilon \left( B + 2B^2A + nB^3A + A(\gamma-1)n\left(B + \frac{3}{2}nB^3\right)\right)$$

$$\leq B, \quad j \in \{1,\dots,n\},$$

$$\| S_n(\rho,u,E) \| \leq \varepsilon \left( B + 2BA + AB^2n + ABn(\gamma-1)\left(B + \frac{3}{2}nB^3\right)\right)$$

$$\leq B.$$

Therefore $S: X \to X$ and it is continuous. Next, $T: X \to Y$ is expansive with a constant $1+\varepsilon$ and if $(\rho^1, u^1, E^1) \in Y$, then

$$T(\rho,u,E) = (\rho^1, u^1, E^1),$$

where

$$(\rho,u,E) = \frac{1}{1+\varepsilon}(\rho^1, u^1, E^1),$$

i.e., $T: X \to Y$ is onto. Hence and Theorem 2.1, we conclude that $T+S$ has a fixed point in $X$. This completes the proof.

**Example 3.1** *Let*

$$\rho_0(x) = E_0(x) = u_{i0}(x) = \frac{1}{1+|x|^2}, \quad i \in \{1,\dots,n\},$$

$x \in \mathbb{R}^n_+$. Here $|x|^2 = x_1^2 + \cdots + x_n^2$. We have

$$|\rho_0(x)|, \quad |E_0(x)|, \quad |u_{i0}(x)| \leq 1, \quad i \in \{1,\dots,n\},$$

$x \in \mathbb{R}^n_+$. Moreover,

$$\left|\frac{\partial}{\partial x_j}\rho_0(x)\right|, \quad \left|\frac{\partial}{\partial x_j}E_0(x)\right|, \quad \left|\frac{\partial}{\partial x_j}u_{i0}(x)\right| = \left|-\frac{2x_j}{(1+|x|^2)^2}\right| \leq 2,$$

$i,j \in \{1,\dots,n\}$, $x \in \mathbb{R}^n_+$. Thus, we can take $B = 2$.

## 4 Conclusions

In this paper we investigate a class of Euler equations. We prove existence of classical solutions for the considered IVP. We propose a new approach which can be used for investigations for existence of classical solutions of other classes of partial differential equations.